\documentclass[a4paper,12pt,reqno]{amsart}
\usepackage{amssymb}
\usepackage{ifthen}
 \usepackage[dvips]{graphicx}
\nonstopmode \numberwithin{equation}{section}
\usepackage{amssymb}
\usepackage{ifthen}
\usepackage{subcaption}
\usepackage{graphicx}
\usepackage{cite}
\usepackage{enumitem}
\usepackage{amsmath}
\usepackage[T1]{fontenc} %skandit
\usepackage[utf8]{inputenc}
\usepackage[usenames,dvipsnames]{color}
\usepackage{color}
\usepackage[english]{babel}
\usepackage{fancyhdr}
\usepackage{fancybox}
\usepackage{tikz}
\nonstopmode\numberwithin{equation}{section}
\newtheorem*{thmA}{Theorem A}
\newtheorem*{thmB}{Theorem B}
\newtheorem*{thmC}{Theorem C}
\newtheorem*{thmD}{Theorem D}

\theoremstyle{plain}

\newtheorem{conj}{Conjecture}

\theoremstyle{definition}
\newtheorem{defn}{Definition}[section]

\newtheorem{thm}{Theorem}[section]
\newtheorem{prob}{Problem}[section]
\newtheorem{cor}{Corollary}[section]
\newtheorem{ques}{Question}[section]
\newtheorem{prop}{Proposition}[section]
\newtheorem{rem}{Remark}[section]
\newtheorem{lem}{Lemma}[section]

\newcounter{minutes}
\divide\time by 60
\newcounter{hours}
\multiply\time by 60
\addtocounter{minutes}{-\time}
\newcounter {own}
\def\theown {\thesection       .\arabic{own}}

\newenvironment{pf}[1][]{%
 \vskip 3mm
 \noindent
 \ifthenelse{\equal{#1}{}}%
  {{\slshape Proof. }}%
  {{\slshape #1.} }%
 }%
{\qed\bigskip}

\newcounter{alphabet}

\def\be{\begin{equation}}
\def\ee{\end{equation}}

\newcommand{\bee}{\begin{enumerate}}
\newcommand{\eee}{\end{enumerate}}

\newcommand{\blem}{\begin{lem}}
\newcommand{\elem}{\end{lem}}
\newcommand{\bthm}{\begin{thm}}
\newcommand{\ethm}{\end{thm}}
\newcommand{\bcor}{\begin{cor}}
\newcommand{\ecor}{\end{cor}}
\newcommand{\beg}{\begin{examp}}
\newcommand{\eeg}{\end{examp}}
\newcommand{\begs}{\begin{examples}}
\newcommand{\eegs}{\end{examples}}

\newcommand{\bdefn}{\begin{defn}}
\newcommand{\edefn}{\end{defn}}

\newcommand{\bprob}{\begin{prob}}
\newcommand{\eprob}{\end{prob}}
\newcommand{\bei}{\begin{itemize}}
\newcommand{\eei}{\end{itemize}}

\newcommand{\bcon}{\begin{conj}}
\newcommand{\econ}{\end{conj}}
\newcommand{\bcons}{\begin{conjs}}
\newcommand{\econs}{\end{conjs}}
\newcommand{\bprop}{\begin{prop}}
\newcommand{\eprop}{\end{prop}}
\newcommand{\br}{\begin{rem}}
\newcommand{\er}{\end{rem}}
\newcommand{\brs}{\begin{rems}}
\newcommand{\ers}{\end{rems}}
\newcommand{\bo}{\begin{obser}}
\newcommand{\eo}{\end{obser}}
\newcommand{\bos}{\begin{obsers}}
\newcommand{\eos}{\end{obsers}}
\newcommand{\bpf}{\begin{pf}}
\newcommand{\epf}{\end{pf}}
\newcommand{\ba}{\begin{array}}
\newcommand{\ea}{\end{array}}
\newcommand{\beq}{\begin{eqnarray}}
\newcommand{\beqq}{\begin{eqnarray*}}
\newcommand{\eeq}{\end{eqnarray}}
\newcommand{\eeqq}{\end{eqnarray*}}

\begin{document}

\title{Bohr Phenomenon for $K$-quasiconformal Harmonic Mappings Involving One Parameter}

\author{Molla Basir Ahamed}
\address{Molla Basir Ahamed, Department of Mathematics, Jadavpur University, Kolkata-700032, West Bengal, India.}
\email{mbahamed.math@jadavpuruniversity.in}

\author{Taimur Rahman}
\address{Taimur Rahman, Department of Mathematics, Jadavpur University, Kolkata-700032, West Bengal, India.}
\email{taimurr.math.rs@jadavpuruniversity.in}

\subjclass[{AMS} Subject Classification:]{Primary 30A10, 30C45, 30C62,  30C50}

\def\thefootnote{}
\footnotetext{ {\tiny File:~\jobname.tex,
printed: \number\year-\number\month-\number\day,
          \thehours.\ifnum\theminutes<10{0}\fi\theminutes }
} \makeatletter\def\thefootnote{\@arabic\c@footnote}\makeatother

\begin{abstract}
In this article, we study Bohr-type inequalities involving a parameter or convex combinations for $K$-quasiconformal, sense-preserving harmonic mappings in $\mathbb{D}$, where the analytic part is subordinate to a convex function. Moreover, we establish similar inequalities when the subordinating function is chosen from the class of concave univalent functions with pole $p$, as well as from the family of concave univalent functions with opening angle $\pi\alpha$. The results generalize several existing results.
\end{abstract}

\maketitle
\pagestyle{myheadings}
\markboth{M. B. Ahamed and T. Rahman}{Bohr's phenomenon for $K$-quasiconformal harmonic mappings}
%\tableofcontents
\section{Introduction}
Bohr’s foundational result from 1914 \cite{Bohr-PLMS-1914}, dealing with power series in complex analysis, gave rise to what is now termed Bohr’s phenomenon. This influential contribution has attracted extensive attention, with significant results established in various function spaces. For comprehensive accounts, we recommend the survey by Abu-Muhanna \emph{et al.} \cite{Abu-Muhanna-Ali-Ponnusami-2016}, the chapter by Garcia \emph{et al.} \cite[Chapter 8]{Garcia-Mashreghi-Ross-2018}, and the references therein.\vspace{2mm}

\subsection{Classical Bohr inequality and its generalizations:}

Let $\mathcal{B}$ represent the class of all analytic functions $f$ on $\mathbb{D}$ bounded by $|f(z)| \leq 1$. The inequality introduced in 1914 by Hardy, after correspondence with Harald Bohr (see \cite{Bohr-PLMS-1914}), is stated as follows:
\begin{thmA}
	If $f(z) = \sum_{n=0}^\infty a_n z^n \in \mathcal{B}$, then
	\begin{align}\label{eq-1.1}
		M_f(r) := \sum_{n=0}^\infty |a_n| r^n \leq 1 \quad \text{for } r \leq \tfrac{1}{3}.
	\end{align}
The radius $1/3$ is best possible.
\end{thmA}

Bohr initially introduced inequality \eqref{eq-1.1} for $r \leq 1/6$. The constant was improved to $r \leq 1/3$ by the independent work of M. Riesz, I. Schur, and F. Wiener, who also confirmed that $1/3$ is the best possible. Consequently, $1/3$ is termed the Bohr radius, while inequality \eqref{eq-1.1} is the Bohr inequality for $\mathcal{B}$. For the function
\begin{align*}
	f_a(z) = \frac{a-z}{1-az} = a - (1-a^2)\sum_{n=1}^{\infty} a^{n-1}z^n, \quad z \in \mathbb{D},
\end{align*}
one easily sees that $M_{f_a}(r) > 1$ holds if and only if $r > 1/(1+2a)$. This implies, as $a \to 1^{-}$, that $1/3$ is indeed sharp.\vspace{2mm}

 It is worth mentioning that in the majorant series $M_f(r)$, if one replaces $|a_0|$ with $|a_0|^2$ in Bohr’s inequality, then the constant $1/3$ improves to $1/2$. Furthermore, when $a_0=0$ in Theorem A, the exact Bohr radius is given by $1/\sqrt{2}$ (see \cite{Kayumov-Ponnusamy-CMFT-2017}, \cite[Corollary 2.9]{Paulsen-Popescu-Singh-PLMS-2002}, and \cite{Ponnusamy-Wirths-CMFT-2020}). Paulsen \emph{et al.} \cite{Paulsen-Singh-2022} later offered an elegant and elementary proof of Bohr’s inequality. In recent years, numerous studies have focused on the Bohr radius and its generalizations. Abu-Muhanna and Ali \cite{Abu-Muhanna-CVEE-2010,Abu-Muhanna-Ali-JMAA-2011} investigated the problem for analytic mappings from $\mathbb{D}$ to simply connected domains and the complement of the unit disk. Further, the Bohr phenomenon in shifted disks and simply connected domains has been studied in \cite{Ahamed-Allu-Halder-AFM-2022,Evdoridis-Ponnusamy-Rasila,Fournier-Ruscheweyh-CRM-2010}. Allu and Halder \cite{Allu-Halder-JMAA-2021} and Bhowmik and Das \cite{Bhowmik-Das-JMAA-2018} studied Bohr-type inequalities in the framework of subordination classes. For a comprehensive overview, we refer to \cite{Abu-Muhanna-CVEE-2010,Abu-Muhanna-Ali-JMAA-2011,Alkhaleefah-Kayumov-Ponnusamy-PAMS-2019,Evdoridis-Ponnusamy-Rasila,Ismagilov-Kayumov-Ponnusamy-JMAA-2020,Kayumov-Ponnusamy-CMFT-2017,Kayumov-Ponnusamy-AASFM-2019} and references therein.
 \subsection{Bohr-Rogosinski inequality and its generalizations:}

 In addition to the Bohr radius, another significant concept is the Rogosinski radius \cite{Rogosinski-MJ-1923}, defined as follows: Let $f(z)=\sum_{n=0}^{\infty}a_nz^n$ be analytic in $\mathbb{D}$ with $|f(z)|\leq 1$ for $z\in\mathbb{D}$. Then, for each $N\geq 1$, the inequality $|S_N(z)|\leq 1$ holds for $|z|\leq 1/2$, where $S_N(z):=\sum_{n=0}^{N-1}a_nz^n$ denotes the $N$-th partial sum. The constant $1/2$ is sharp. Motivated by this notion, Kayumov and Ponnusamy \cite{Kayumov-Ponnusamy} introduced the Bohr-Rogosinski sum
 \begin{align*}
 	R^{f}_N(z):=|f(z)|+\sum_{n=N}^{\infty}|a_n||z|^n,
 \end{align*}
 which clearly satisfies $|S_N(z)|=|f(z)-\sum_{n=N}^{\infty}a_nz^n|\leq R^{f}_N(z)$. Furthermore, this sum generalizes the classical Bohr sum, corresponding to the case $N=1$ with $f(0)$ replaced by $f(z)$. The Bohr-Rogosinski radius is then defined as the maximal $r\in(0,1)$ such that $R^{f}_N(z)\leq 1$ for $|z|\leq r$.
 
\begin{thmB}\cite{Kayumov-Ponnusamy}
	Let $f(z)=\sum_{n=0}^{\infty}a_nz^n$ be an analytic function in $\mathbb{D}$ and $|f(z)|\leq 1$. Then, for each $N\in\mathbb{N}$, we have 
	\begin{align}
		|f(z)|+\sum_{n=N}^{\infty}|a_n||z|^n\leq 1
	\end{align}
	for $|z|=r\leq R_N$, where $R_N$ is the positive root of the equation $2(1+r)r^N-(1-r)^2=0$. The radius $R_N$ is the best possible.\vspace{2mm}
\end{thmB}
Recently, Hamada \emph{et al.} \cite{Hamada-Honda-Kohr-AMP-2025} studied the Bohr--Rogosinski radius for holomorphic mappings defined on the unit ball of a complex Banach space, taking values in higher dimensional Banach spaces. Furthermore, they established the Bohr--Rogosinski radius for certain subordination classes on the unit ball of a complex Banach space. For recent study with various class of functions on Bohr-Rogosinski inequality, we refer to the articles \cite{Ahamed-CMFT-2022,Gangania-Kumar-MJM-2022,Allu-Arora-JMAA-2023,Kayumov-Ponnusamy} \vspace{2mm}

We recall the notion of differential subordination, which has proven to be an effective approach in solving numerous problems in geometric function theory. A detailed discussion of subordination classes is found in \cite[Chapter 6]{Duren-1983}.
\begin{defn}
Suppose that $f$ and $g$ are analytic in $\mathbb{D}$. Then $f$ is said to be subordinate to $g$, expressed as $f\prec g$, if there exists a Schwarz function $w\in\mathcal{B}$ satisfying $w(0)=0$ such that $f(z)=g(w(z))$ for all $z\in\mathbb{D}$. If $g$ is univalent in $\mathbb{D}$, this subordination holds exactly when $f(0)=g(0)$ and $f(\mathbb{D})\subset g(\mathbb{D})$. Moreover, the condition $f\prec g$ guarantees that $|f^{\prime}(0)|\leq |g^{\prime}(0)|$. 
\end{defn}
 
\subsection{Distance formulation of the Bohr inequality:}

The investigation of the Bohr phenomenon for subordinated functions was first carried out by Abu-Muhanna \cite{Abu-Muhanna-CVEE-2010}. In this context, we denote by $S(g):=\{f:f\prec g\}$ the family of all analytic functions $f$ in $\mathbb{D}$ that are subordinate to a univalent function $g$. The class $S(g)$ is said to satisfy the Bohr phenomenon if there exists $r_g\in(0,1)$ such that, for every $f(z)=\sum_{n=0}^{\infty}a_nz^n\in S(g)$, one has 
\begin{align}\label{eq-1.3}
	\sum_{n=1}^{\infty}|a_n|r^n\leq d(g(0),\partial g(\mathbb{D}))\quad \text{for all } |z|=r\leq r_g,
\end{align}
where $d(g(0),\partial g(\mathbb{D}))$ represents the Euclidean distance between $g(0)$ and $\partial g(\mathbb{D})$.\vspace{2mm}

The Bohr radius for the class $S(g)$ is defined to be the largest $r_g$ such that the Bohr inequality holds. Similarly, $S(g)$ is said to exhibit the Bohr--Rogosinski property if there exists some $r_{N,g}\in(0,1]$ with the property that
\begin{align}\label{eq-1.4}
	|f(z)|+\sum_{n=N}^{\infty}|a_n|r^n \leq |g(0)|+d(g(0),\partial g(\mathbb{D})),
\end{align}
whenever $|z|=r\leq r_{N,g}$ for all $f(z)=\sum_{n=0}^{\infty}a_nz^n\in S(g)$. The largest such $r_{N,g}$ is called the Bohr--Rogosinski radius.\vspace{2mm}

For the mapping $g(z)=(a-z)/(1-\bar{a}z)$ with $|a|\leq 1$, we see that $g(\mathbb{D})=\mathbb{D}$, hence $S(g)=\mathcal{B}$ and $d(g(0),\partial g(\mathbb{D}))=1-|a|$. Therefore, from \eqref{eq-1.1}, inequality \eqref{eq-1.3} holds for $|z|\leq 1/3$, while Theorem A shows that \eqref{eq-1.4} is true with $r_{N,g}=R_N$.

\subsection{Harmonic Mappings:}
For a simply connected domain $\Omega$, let $f=u+iv$ be a complex-valued function. 
The function is harmonic in $\Omega$ if and only if it satisfies $4f_{z\bar{z}}=0$. 
This implies that $u$ and $v$ are harmonic functions. 
It is well known that a harmonic function admits a decomposition $f(z)=h(z)+\overline{g(z)}$, where $h$ and $g$ are analytic. 
Here $h$ is called the analytic part and $g$ the co-analytic part of $f$. 
This splitting is unique up to an additive constant \cite{Duren-2004}. 
A classical result due to Lewy \cite{Lewy-BAMS-1936} shows that $f$ is locally univalent if and only if $J_f(z)=|h'(z)|^2-|g'(z)|^2\neq 0$.
Moreover, $f$ is locally univalent and preserves orientation when $J_f(z)>0$, or equivalently when $h'(z)\neq 0$ and $|\omega|<1$, with $\omega=g'/h'$.

\subsection{K-quasiconformal Harmonic Mappings:}

A harmonic function $f=h+\bar{g}$ is called $K$-quasiconformal harmonic in $\mathbb{D}$ if it is locally univalent, sense-preserving, and satisfies $|\omega_f(z)|\leq k<1$ for all $z\in\mathbb{D}$, where $K=(1+k)/(1-k)\geq 1$ (see \cite{Kalaj-MZ-2008,Martio-AASFAI-1968}). Clearly, as $k\to 1$, it follows that $K\to\infty$. In 2018, Ponnusamy \emph{et al.}\,\cite{Kayumov-Ponnusamy-Shakirov-MN-2018} studied the class of sense-preserving harmonic mappings $f=h+\bar{g}$ in the unit disk $\mathbb{D}$, where $h$ and $g$ are analytic with $g(0)=0$, and determined the Bohr radius under certain assumptions on $h$ and $g$. Later, in 2019, Liu and Ponnusamy \cite{Liu-Ponnusamy-BMMSS-2019} introduced the Bohr radius for $K$-quasiconformal, sense-preserving harmonic functions in $\mathbb{D}$, where $h$ is subordinate to some analytic function $\phi$, and proposed two conjectures. Subsequently, Liu \emph{et al.}\,\cite{Liu-Ponnusami-Wang-2020} made a significant contribution by investigating the Bohr radius for sense-preserving $K$-quasiregular harmonic mappings in $\mathbb{D}$, where $h(z)-h(0)$ is quasi-subordinate to a given analytic function. They not only provided sharper versions of four theorems from \cite{Liu-Ponnusamy-BMMSS-2019}, but also resolved the two conjectures posed therein.\vspace{2mm}

In 2019, Liu and Ponnusamy \cite{Liu-Ponnusamy-BMMSS-2019} established results concerning the Bohr radii of harmonic mappings in the unit disk $\mathbb{D}$, where the analytic part is subordinate to a given analytic function.

\begin{thmC}\cite{Liu-Ponnusamy-BMMSS-2019}
	Suppose that $f(z)=h(z)+\overline{g(z)}=\sum_{n=0}^{\infty}a_nz^n+\overline{\sum_{n=1}^{\infty}b_nz^n}$ is a sense-preserving $K$-quasiconformal harmonic mapping in $\mathbb{D}$ and $h\prec\phi$, where $\phi$ is univalent and convex in $\mathbb{D}$. Then
	\begin{align*}
		\sum_{n=1}^{\infty}(|a_n|+|b_n|)r^n\leq d(\phi(0),\partial\phi(\mathbb{D}))\;\; \mbox{for}\;\; r\leq \frac{K+1}{5K+1},
	\end{align*}
	The number $(K+1)/(5K+1)$ is sharp.
\end{thmC}
\begin{thmD}\cite{Liu-Ponnusamy-BMMSS-2019}
	Suppose that $f(z)=h(z)+\overline{g(z)}=\sum_{n=0}^{\infty}a_nz^n+\overline{\sum_{n=1}^{\infty}b_nz^n}$ is a sense-preserving $K$-quasiconformal harmonic mapping in $\mathbb{D}$ and $h\prec\phi$, where $\phi$ is analytic and univalent in $\mathbb{D}$. Then
	\begin{align*}
		\sum_{n=1}^{\infty}(|a_n|+|b_n|)r^n\leq d(\phi(0),\partial\phi(\mathbb{D}))\;\; \mbox{for}\;\; r\leq r_u(k),
	\end{align*}
	where $k=(K-1)/(K+1)$ and $r_u(k)\in(0,1)$ is the positive root of the equation 
	\begin{align*}
		(1-r)^2-4r(1+k\sqrt{1+r})=0.
	\end{align*}
\end{thmD}
In \cite{Long-Wang-Wu-2022}, Long \emph{et al.} established Bohr-type inequalities with a parameter or involving convex combinations for bounded analytic functions. Motivated by these results, Arora and Vinayak \cite{Arora-Vinayak-CVVE-2025} investigated Bohr-type inequalities for several classes of analytic functions. They began by studying the Bohr phenomenon through convex combinations and extended their work to other Bohr-type inequalities for univalent analytic functions. Furthermore, they examined similar inequalities for the Ma-Minda classes of convex and starlike functions, and also provided the same estimates for concave univalent functions.\vspace{2mm}

The results in \cite{Long-Wang-Wu-2022, Arora-Vinayak-CVVE-2025} motivate us to pose the following questions.
\begin{ques}\label{ques-1.1}
	Is it possible to establish Bohr-type inequalities involving parameters or convex combinations for $K$-quasiconformal, sense-preserving harmonic functions in $\mathbb{D}$, assuming the analytic part is subordinate to a convex function?
\end{ques}

\begin{ques}\label{ques-1.2}
	Is it possible to obtain Bohr-type inequalities with a parameter or convex combination for $K$-quasiconformal harmonic mappings whose analytic part is subordinate to a concave univalent function with pole $p$?
\end{ques}
\begin{ques}\label{ques-1.3}
	Is it possible to derive Bohr-type inequalities with a parameter or convex combination for $K$-quasiconformal harmonic mappings whose analytic component is subordinate to a concave univalent function of opening angle $\pi\alpha$?
\end{ques}
As part of our main contributions, we provide affirmative answers to Questions \ref*{ques-1.1}, \ref*{ques-1.2}, and \ref*{ques-1.3}, which are discussed in Section2, Section3, and Section~4, respectively.\vspace{2mm}

We organize the paper as follows: In Section2, we study Bohr-type inequalities involving parameters or convex combinations for $K$-quasiconformal harmonic mappings subordinate to convex functions. In Section 3, we derive similar inequalities for the case when the subordinate function belongs to the family of concave univalent functions with pole $p$. In Section 4, We provide Bohr-type inequalities assuming the subordinate function belongs to the family of concave univalent functions with opening angle $\pi\alpha$.
\section{\bf{Bohr's phenomenon related to the family of convex functions}}
Convex functions form an important class of functions in geometric function
theory. A single-valued function $f$ is said to be univalent in a domain 
$D \subset \mathbb{C}$ if, for any two distinct points $z_1, z_2 \in D$, their 
corresponding function values $f(z_1)$ and $f(z_2)$ are also distinct. A univalent 
function $f$ is said to be convex in a domain $D \subset \mathbb{C}$ if it is analytic 
and the image $f(D)$ is a convex set. A set in the complex plane is convex if, for any 
two points $w_1, w_2 \in f(D)$, the line segment joining them lies entirely within $f(D)$. 
Equivalently, a univalent function $f$ is convex in $D$ if it maps $D$ conformally onto a convex domain. It is worth noting that the function defined by 
\begin{align*}
	c(z):=\frac{z}{1-z}, \;\; z\in \mathbb{D},
\end{align*} 
is a classical example of a convex function and often plays the role of an extremal function.\vspace{2mm}

Before presenting our main results, we recall the following lemmas.
\begin{lem}\label{lem-2.1}\cite{Gangania-Kumar-MJM-2022}
	Suppose that $h(z)=\sum_{n=0}^{\infty}a_nz^n$ and $g(z)=\sum_{n=0}^{\infty}b_nz^n$ are two analytic functions in $\mathbb{D}$ and $h\prec g$, then for $N\in\mathbb{N}$, we have
	\begin{align*}
		\sum_{n=N}^{\infty}|a_n|r^n\leq\sum_{n=N}^{\infty}|b_n|r^n \;\;\; \mbox{holds for}\;\; |z|=r\leq 1/3.
	\end{align*}
\end{lem}
\begin{lem}\label{lem-2.2}\cite{Alkhaleefah-Kayumov-Ponnusamy-PAMS-2019,Liu-Ponnusami-Wang-2020}
	Suppose that $h(z)=\sum_{n=0}^{\infty}a_nz^n$ and $g(z)=\sum_{n=0}^{\infty}b_nz^n$ are two analytic functions in $\mathbb{D}$. If $|g^{\prime}(z)|\leq k|h^{\prime}(z)|$ in $\mathbb{D}$ for some $k\in(0,1]$, then
	\begin{align*}
		\sum_{n=1}^{\infty}|b_n|r^n\leq k\sum_{n=1}^{\infty}|a_n|r^n \;\; \mbox{for}\;\; |z|=r\leq 1/3.
	\end{align*}
\end{lem}

\begin{lem}\cite{Abu-Muhanna-CVEE-2010}\label{lem-2.3}
	Suppose $g(z)=\sum_{n=0}^{\infty}a_nz^n$ is an analytic and univalent map from $\mathbb{D}$ onto a convex domain $g(\mathbb{D})$, then 
	\begin{align*}
		\frac{1}{2}(1-|z|^2)|g^{\prime}(z)|\leq d(g(z),\partial g(\mathbb{D}))\leq (1-|z|^2)|g^{\prime}(z)|\;\;\; \mbox{for}\;\; z\in\mathbb{D}.
	\end{align*}
	If $f(z)=\sum_{n=0}^{\infty}b_nz^n\prec g(z)$, then 
	\begin{align*}
		|b_n|\leq |g^{\prime}(0)|\leq 2d(g(0),\partial g(\mathbb{D}))\;\;\; \mbox{for}\;\; n\geq 1.
	\end{align*}
\end{lem}
We now present our first result on Bohr’s phenomenon, formulated in terms of convex combinations for convex functions.

\begin{thm}\label{thm-2.1}
	Suppose that $f(z)=h(z)+\overline{g(z)}=\sum_{n=0}^{\infty}a_nz^n+\overline{\sum_{n=1}^{\infty}b_nz^n}$ is a sense-preserving $K$-quasiconformal harmonic mapping in $\mathbb{D}$ and $h\prec\phi$, where $\phi$ is convex function. Then for arbitrary $t\in[0,1]$ and $m\in\mathbb{N}$, we have
	\begin{align*}
		t|h(z^m)|+(1-t)\sum_{n=0}^{\infty}(|a_n|+|b_n|)r^n\leq|\phi(0)|+ d(\phi(0),\partial\phi(\mathbb{D})) 
	\end{align*}
for $r\leq \min\{1/3,R_{K,m,t}\}$, where $R_{K,m,t}$ is the unique root in $(0,1)$ of the equation 
	\begin{align}\label{eq-2.1}
		D^t_{K,m}(r):=\frac{tr^m}{1-r^m}+\frac{2(1-t)Kr}{(K+1)(1-r)}-\frac{1}{2}=0.
	\end{align}
	The radius $R_{K,m,t}$ is sharp if $R_{K,m,t}\leq 1/3$.
\end{thm}
The following corollary is obtained from the above theorem by taking $K=1$.
\begin{cor}
 Let $f(z)=\sum_{n=0}^{\infty}a_nz^n$ be a an analytic function defined on $\mathbb{D}$ and $f\prec\phi$, where $\phi$ is a convex function. Then for arbitrary $t\in[0,1]$ and $m\in\mathbb{N}$, we have 
 \begin{align*}
 	t|f(z^m)|+(1-t)\sum_{n=0}^{\infty}|a_n|r^n\leq |\phi(0)|+d(\phi(0),\partial\phi(\mathbb{D}))
 \end{align*}
 for $r\leq \min\{1/3,R_{m,t}\}$, where $R_{m,t}$ is the unique in $(0,1)$ of the equation 
 \begin{align*}
 	D^t_m(r):=\frac{tr^m}{1-r^m}+\frac{(1-t)r}{1-r}-\frac{1}{2}=0.
 \end{align*}
 The number $R_{m,t}$ cannot be improved if $R_{m,t}\leq1/3$.
\end{cor}
\begin{rem}
	If we set $t=0$ in Theorem \ref{thm-2.1}, it coincides with Theorem C.
\end{rem}

\begin{proof}[\bf Proof of Theorem \ref{thm-2.1}]
	Since $h\prec\phi$ and $\phi$ is analytic convex univalent function in $\mathbb{D}$, it follows from Lemma \ref{lem-2.3} that $|a_n|\leq |\phi^{\prime}(0)|\leq 2d(\phi(0),\partial \phi(\mathbb{D}))$ for $n\geq 1$. Then, we have 
	\begin{align}\label{eq-2.2}
		\sum_{n=1}^{\infty}|a_n|r^n\leq 2d(\phi(0),\partial \phi(\mathbb{D}))\sum_{n=1}^{\infty}r^n=2d(\phi(0),\partial \phi(\mathbb{D}))\frac{r}{1-r}.
	\end{align}
	Since $f$ is $K$-quasiconformal sense-preserving harmonic mapping on $\mathbb{D}$, it follows from Schwarz’s Lemma that the dilatation $\omega = g'/h'$ is analytic with $|\omega(z)| \leq k$ in $\mathbb{D}$, where $K = (1 + k)/(1 - k)$ and $k \in [0,1)$. Hence, by using Lemma \ref{lem-2.2}, we obtain
	\begin{align}\label{eq-2.3}
		\sum_{n=1}^{\infty}|b_n|r^n\leq k\sum_{n=1}^{\infty}|a_n|r^n\leq 2kd(\phi(0),\partial \phi(\mathbb{D}))\frac{r}{1-r}.
	\end{align}
	Since $\phi$ is univalent and $h\prec\phi$, it follows that $h(0)=\phi(0)$ and hence
	\begin{align}\label{eq-2.4}
		|h(z^m)|\leq |h(0)|+\left|\sum_{n=1}^{\infty}a_nz^{mn}\right|\leq |\phi(0)|+2d(\phi(0),\partial \phi(\mathbb{D}))\frac{r^m}{1-r^m}.
	\end{align}
	From inequalities \eqref{eq-2.2}, \eqref{eq-2.3}, and \eqref{eq-2.4}, it follows that 
	\begin{align*}
		t|h(z^m)|+&(1-t)\sum_{n=0}^{\infty}(|a_n|+|b_n|)r^n\\\leq& |\phi(0)|+2d(\phi(0),\partial \phi(\mathbb{D}))\left(\frac{tr^m}{1-r^m}+\frac{2(1-t)Kr}{(K+1)(1-r)}\right)\\=&|\phi(0)|+2d(\phi(0),\partial \phi(\mathbb{D}))\left(D^t_{K,m}(r)+\frac{1}{2}\right),
	\end{align*}
	where $D^t_{K,m}(r)$ is defined by \eqref{eq-2.1}. It is easy to observe that the function $D^t_{K,m}(r)$ increases strictly on the interval $(0,1)$, with $D^t_{K,m}(0)=-1/2<0$ and $\lim_{r\to1} D^t_{K,m}(r)=+\infty$. It follows that the equation $D^t_{K,m}(r)=0$ admits exactly one solution $R_{K,m,t}$ in $(0,1)$. Consequently, $D^t_{K,m}(r)\leq 0$ if, and only if, $r\leq\min\{1/3, R_{K,m,t}\}$, proving the desired inequality.\vspace{2mm}
	
	For the proof of the sharpness in the case $R_{K,m,t}\leq 1/3$, we consider the function $f(z)=h(z)+\overline{g(z)}$ in $\mathbb{D}$ such that 
	\begin{align*}
		\phi(z)=h(z)=\frac{z}{1-z}=\sum_{n=1}^{\infty}z^n \;\;\;\mbox{and}\;\; g(z)=k\mu\sum_{n=1}^{\infty}z^n,
	\end{align*}
	where $|\mu|=1$ and $k=(K-1)/(K+1)$. It is well-known that $d(\phi(0),\partial \phi(\mathbb{D}))=1/2$ and $|\phi(0)|=0$. For $z=r$, we have
	\begin{align*}
		t|h(r^m)|+(1-t)\sum_{n=0}^{\infty}(|a_n|+|k\mu a_n|)r^n&=\frac{tr^m}{1-r^m}+\frac{2(1-t)Kr}{(K+1)(1-r)}\\&=D^t_{K,m}(r)+\frac{1}{2}\\&=D^t_{K,m}(r)+|\phi(0)|+d(\phi(0),\partial \phi(\mathbb{D}))\\&>|\phi(0)|+d(\phi(0),\partial \phi(\mathbb{D}))
	\end{align*}
	holds true for $r>R_{K,m,t}$. This completes the proof.
\end{proof}
Table~1 provides the values of $R_{K,m,t}$ for selected values of $K \geq 1$, $m \in \mathbb{N}$, and $t \in [0,1]$, and Figure~1 indicates their corresponding locations.
\begin{table}[ht]
	
	\centering
	\begin{tabular}{|l|l|l|l|l|l|l|l|l|l|}
		\hline
		$K$& $m$&$t$& $R_{K,m,t}$ \\
		\hline
		$1$& $1$&$0.2$& $0.2941$ \\
		\hline
		$3$& $2$&$0.5$& $1/3$ \\
		\hline
		$9$& $5$&$0.8$& $0.2573$ \\
		\hline
		$1$& $1$&$1 $& $0.2$ \\
		\hline
		$25$& $15$& $2$& $0.1150$\\
		\hline	
	\end{tabular}
	\caption{The unique positive root of Equation \eqref{eq-2.1} in $(0,1)$, denoted as $R_{K,m,t}$, is shown in the table.}
\end{table}

For the same class, we now examine a Bohr-type inequality involving a parameter, which leads to the following result.
 
\begin{thm}\label{thm-2.2}
	Suppose that $f(z)=h(z)+\overline{g(z)}=\sum_{n=0}^{\infty}a_nz^n+\overline{\sum_{n=1}^{\infty}b_nz^n}$ is a sense-preserving $K$-quasiconformal harmonic mapping in $\mathbb{D}$ and $h\prec\phi$,  where $\phi$ is analytic convex univalent function. Then for arbitrary $\lambda>0$ and $m\in\mathbb{N}$, the inequality
	\begin{align*}
		|h(z^m)|+\lambda\sum_{n=1}^{\infty}(|a_n|+|b_n|)r^n\leq|\phi(0)|+ d(\phi(0),\partial\phi(\mathbb{D}))
	\end{align*}
	holds for $|z|=r\leq R^*_{K,m,\lambda}$, which is the unique root in $(0,1)$ of the equation 
	\begin{align}\label{eq-2.5}
		E^{\lambda}_{K,m}(r):=\frac{r^m}{1-r^m}+\frac{2 K\lambda r}{(K+1)(1-r)}-\frac{1}{2}=0.
	\end{align}
	The radius $R^*_{K,m,\lambda}$ is sharp if $R^*_{K,m,\lambda}\leq1/3$.
\end{thm}
If we set $K=1$, the above theorem reduces to the following corollary.
\begin{cor}\label{cor-2.2}
Let $f(z)=\sum_{n=0}^{\infty}a_nz^n$ be a an analytic function defined on $\mathbb{D}$ and $f\prec\phi$, where $\phi$ is a convex function. Then for arbitrary $\lambda>0$ and $m\in\mathbb{N}$, we have 
\begin{align*}
	|f(z^m)|+\lambda\sum_{n=1}^{\infty}|a_n|r^n\leq |\phi(0)|+d(\phi(0),\partial\phi(\mathbb{D}))
\end{align*}
for $r\leq \min\{1/3,R_{m,\lambda}\}$, where $R^*_{m,\lambda}$ is the unique in $(0,1)$ of the equation 
\begin{align*}
	E^{\lambda}_m(r):=\frac{r^m}{1-r^m}+\frac{\lambda r}{1-r}-\frac{1}{2}=0.
\end{align*}
The number $R^*_{m,\lambda}$ cannot be improved if $R^*_{m,\lambda}\leq 1/3$.
\end{cor}
\begin{rem} 

\begin{enumerate}
	\item[(i)]  Under the assumptions $m \to \infty$ and $\lambda=1$ in Corollary \ref{cor-2.2}, the radius reduces to $r=1/3$, coinciding with the conclusion of \cite[Remark 1]{Abu-Muhanna-CVEE-2010}.\vspace{1.2mm}
	
	\item[(ii)] If we set $m=\lambda=1$ in Corollary \ref{cor-2.2}, the radius reduces to $r=1/5$, which coincides with the result in \cite[Corollary 3]{Kayumov-Khammatova-Ponnusamy-JMAA-2021} for $N=1$.
\end{enumerate}
\end{rem}
\begin{proof}[\bf Proof of Theorem \ref{thm-2.2}]
	By applying the same method as in the proof of Theorem \ref{thm-2.1}, and using inequalities \eqref{eq-2.2}, \eqref{eq-2.3}, and \eqref{eq-2.4}, we find that
	\begin{align*}
		|h(z^m)|+\lambda\sum_{n=1}^{\infty}(|a_n|+|b_n|)r^n\leq&|\phi(0)|+2d(\phi(0),\partial \phi(\mathbb{D}))\left(\frac{r^m}{1-r^m}+\frac{2 K\lambda r}{(K+1)(1-r)}\right)\\=&|\phi(0)|+2d(\phi(0),\partial \phi(\mathbb{D}))\left(E^{\lambda}_{K,m}(r)+\frac{1}{2}\right),
	\end{align*}
	where $E^{\lambda}_{K,m}(r)$ is defined by \eqref{eq-2.5}. The function $E^{\lambda}_{K,m}(r)$ is clearly increasing on the interval $(0,1)$, with $E^{\lambda}_{K,m}(0)<0$ and $\lim_{r\to 1}E^{\lambda}_{K,m}(r)=+\infty$. Therefore, the equation $E^{\lambda}_{K,m}(r)=0$ admits a unique solution in $(0,1)$, which we denote by $R^*_{K,m,\lambda}$. Consequently, the last expression is less than $|\phi(0)|+d(\phi(0),\partial\phi(\mathbb{D}))$ for all $r\leq \min\{1/3,R^*_{K,m,\lambda}\}$.\vspace{1.2mm}

	The sharpness can be established by proceeding exactly as in the sharpness proof of Theorem \ref{thm-2.1}. Thus, we omit the repetition of those details.

\end{proof}
Table 2 displays the values of $R^*_{K,m,\lambda}$ for chosen values of $K \geq 1$, $m \in \mathbb{N}$, and $\lambda > 0$, with their corresponding locations shown in Figure 2.
\begin{table}[ht]
	\centering
	\begin{tabular}{|l|l|l|l|l|l|l|l|l|l|}
		\hline
		$K$& $m$&$\lambda$& $R^*_{K,m,\lambda}$ \\
		\hline
		$1$& $1$&$0.2$& $0.2941$ \\
		\hline
		$3$& $2$&$0.5$& $1/3$ \\
		\hline
		$9$& $5$&$0.8$& $0.2573$ \\
		\hline
		$1$& $1$&$1 $& $0.2$ \\
		\hline
		$25$& $15$& $2$& $0.1150$\\
		\hline
	\end{tabular}
	\caption{In this table, it is shown that $R^*_{K,m,\lambda}$ is the unique positive root of the equation \eqref{eq-2.5} in $(0,1)$.}
\end{table}

Replacing the coefficient $a_1$ with $h^{\prime}(z^m)$ in the Taylor series and examining a Bohr-type inequality involving a parameter, we obtain the following result.

\begin{thm}\label{thm-2.3}
	Suppose that $f(z)=h(z)+\overline{g(z)}=\sum_{n=0}^{\infty}a_nz^n+\overline{\sum_{n=1}^{\infty}b_nz^n}$ is a sense-preserving $K$-quasiconformal harmonic mapping in $\mathbb{D}$ and $h\prec\phi$,  where $\phi$ is analytic convex univalent function. Then for arbitrary $\lambda>0$ and $m\in\mathbb{N}$, the inequality
	\begin{align*}
		|h(z^m)|+|h^{\prime}(z^m)|r+\lambda\left(\sum_{n=2}^{\infty}|a_n|r^n+\sum_{n=1}^{\infty}|b_n|r^n\right)\leq|\phi(0)|+ d(\phi(0),\partial\phi(\mathbb{D}))
	\end{align*}
	holds for $|z|=r\leq \min\{1/3,R^{**}_{K,m,\lambda}\}$, which is the unique root in $(0,1)$ of the equation 
	\begin{align}\label{eq-2.6}
		F^{\lambda}_{K,m}(r):=\frac{r^m}{1-r^m}+\frac{r}{(1-r^m)^2}+\frac{\lambda(K-1+r(K+1))r}{(K+1)(1-r)}-\frac{1}{2}=0.
	\end{align}
	The number $R^{**}_{K,m,\lambda}$ is sharp if $R^{**}_{K,m,\lambda}\leq 1/3$.
\end{thm}
The corollary that follows arises as a special case of the above theorem with $K=1$.
\begin{cor}\label{cor-2.3}
	Let $f(z)=\sum_{n=0}^{\infty}a_nz^n$ be a an analytic function defined on $\mathbb{D}$ and $f\prec\phi$, where $\phi$ is a convex function. Then for arbitrary $\lambda>0$ and $m\in\mathbb{N}$, the inequality
	\begin{align*}
		|f(z^m)|+|f^{\prime}(z^m)|r+\lambda\sum_{n=2}^{\infty}|a_n|r^n\leq|\phi(0)|+ d(\phi(0),\partial\phi(\mathbb{D}))
	\end{align*}
	holds for $|z|=r\leq \min\{1/3,R^{**}_{m,\lambda}\}$, which is the unique root in $(0,1)$ of the equation 
	\begin{align*}
		F^{\lambda}_{m}(r):=\frac{r^m}{1-r^m}+\frac{r}{(1-r^m)^2}+\frac{\lambda r^2}{1- r}-\frac{1}{2}=0.
	\end{align*}
	The number $R^{**}_{m,\lambda}$ is sharp if $R^{**}_{m,\lambda}\leq1/3$.
\end{cor}
\begin{rem}
	By setting $m \to \infty$ and $\lambda = 1$ in Corollary \ref{cor-2.3}, we find the radius $r = 1/3$, which matches the conclusion of \cite[Remark 1]{Abu-Muhanna-CVEE-2010}.
\end{rem}
\begin{proof}[\bf Proof of Theorem \ref{thm-2.3}]
	Given that $h\prec\phi$ and $\phi$ is analytic convex univalent function in $\mathbb{D}$, Lemma \ref{lem-2.3} implies $|a_n|\leq |\phi^{\prime}(0)|\leq 2d(\phi(0),\partial \phi(\mathbb{D}))$ for $n\geq 1$. Consequently, we obtain 
	\begin{align}\label{eq-2.7}
		\sum_{n=2}^{\infty}|a_n|r^n\leq 2d(\phi(0),\partial \phi(\mathbb{D}))\sum_{n=2}^{\infty}r^n=2d(\phi(0),\partial \phi(\mathbb{D}))\frac{r^2}{1-r}.
	\end{align}
	Following the same approach as in Theorem \ref{thm-2.1}, and making use of Lemma \ref{lem-2.2}, we obtain the inequalities \eqref{eq-2.3} and \eqref{eq-2.4}.\vspace{1.2mm}
	
	A simple calculation shows that
	\begin{align}\label{eq-2.8}
		|h^{\prime}(z^m)|=\left|\sum_{n=1}^{\infty}na_nz^{m(n-1)}\right|&\leq 2d(\phi(0),\partial \phi(\mathbb{D}))\sum_{n=1}^{\infty}nr^{m(n-1)}\\&= 2d(\phi(0),\partial \phi(\mathbb{D}))\frac{1}{(1-r^m)^2}\nonumber.
	\end{align}
	Applying the inequalities \eqref{eq-2.3}, \eqref{eq-2.4}, \eqref{eq-2.7} and \eqref{eq-2.8}, we obtain
	\begin{align*}
		|h&(z^m)|+|h^{\prime}(z^m)|r+\lambda\left(\sum_{n=2}^{\infty}|a_n|r^n+\sum_{n=1}^{\infty}|b_n|r^n\right)\\\leq&|\phi(0)|+2d(\phi(0),\partial \phi(\mathbb{D}))\left(\frac{r^m}{1-r^m}+\frac{r}{(1-r^m)^2}+\frac{\lambda(k+r)r}{1-r}\right)\\=&|\phi(0)|+2d(\phi(0),\partial \phi(\mathbb{D}))\left(F^{\lambda}_{K,m}(r)+\frac{1}{2}\right),
	\end{align*} 
	where the function $F^{\lambda}_{K,m}(r)$ is defined in \eqref{eq-2.6}. It is easy to observe that the function $F^{\lambda}_{K,m}(r)$  increases on the interval $(0,1)$, with $F^{\lambda}_{K,m}(0)=-1/2<0$ and $\lim_{r\to1}F^{\lambda}_{K,m}(r)=+\infty$. Thus, the equation $F^{\lambda}_{K,m}(r)=0$ has a unique root in $(0,1)$, denoted by $R^{**}_{K,m,\lambda}$. Consequently, we find that $F^{\lambda}_{K,m}(r)\leq 0$ holds if and only if $r\leq \min\{1/3,R^{**}_{K,m,\lambda}\}$ and this gives the desired inequality.\vspace{2mm}
	
	To prove the sharpness part of the result in the case $R^{**}_{K,m,\lambda}\leq 1/3$, we consider the same function used in the sharpness portion of the proof of Theorem \ref{thm-2.1}. For $z=r$, we obtain
	\begin{align*}
		|&h(r^m)|+|h^{\prime}(r^m)|r+\lambda\left(\sum_{n=2}^{\infty}|a_n|r^n+\sum_{n=1}^{\infty}|\mu k a_n|r^n\right)\\&=\frac{r^m}{1-r^m}+\frac{r}{(1-r^m)^2}+\frac{\lambda(k+r)r}{1-r}\\&=F^{\lambda}_{K,m}(r)+\frac{1}{2}\\&=F^{\lambda}_{K,m}(r)+|\phi(0)|+d(\phi(0),\partial \phi(\mathbb{D}))\\&>|\phi(0)|+d(\phi(0),\partial \phi(\mathbb{D}))
	\end{align*}
	The last inequality holds only when $r>R^{**}_{K,m,\lambda}$. This completes the proof.
	
\end{proof}
The values of $R^{**}_{K,m,\lambda}$ for certain $K \geq 1$, $m \in \mathbb{N}$, and $\lambda > 0$ are listed in Table 2, and their locations appear in Figure 2.
\begin{table}[ht]
	\centering
	\begin{tabular}{|l|l|l|l|l|l|l|l|l|l|}
		\hline
		$K$& $m$&$\lambda$& $R^{**}_{K,m,\lambda}$ \\
		\hline
		$6$& $3$&$0.2$& $0.3321$ \\
		\hline
		$9$& $5$&$0.5$& $0.2832$ \\
		\hline
		$25$& $15$&$1$&  $0.2063$ \\
		\hline
		$20$& $5$&$2 $& $0.1447$ \\
		\hline
		$15$& $10$& $5$ &$0.0807$\\
		\hline
		
	\end{tabular}
	\caption{The table indicates that for various values of $K$, $m$, and $\lambda$, $R^{**}_{K,m,\lambda}$ is the unique positive root of \eqref{eq-2.6} on the interval $(0,1)$.}
	\end{table}
\section{\bf{Bohr-type inequalities for the family of concave univalent functions with pole $p$}}
5. Let $\widehat{\mathbb{C}} := \mathbb{C} \cup \{\infty\}$ denote the extended complex plane. We define $\widehat{C_p}$ to be the family of meromorphic univalent mappings $f : \mathbb{D} \to \widehat{\mathbb{C}}$ such that:
\begin{enumerate}
	\item[(i)] $f$ is analytic in $\mathbb{D} \setminus \{p\}$, and the complement $\widehat{\mathbb{C}} \setminus f(\mathbb{D})$ is convex.\vspace{1.3mm}
	\item[(ii)] $f$ has a simple pole at $p$.
\end{enumerate}

Without loss of generality, a suitable rotation allows us to suppose that $0 < p < 1$. Each $f \in \widehat{C_p}$ is a concave univalent mapping with a pole at $p \in (0,1)$, possessing the Taylor expansion
 \[
 f(z) = \sum_{n=0}^{\infty} a_n z^n \quad z \in \mathbb{D}_p,
 \]
 where $\mathbb{D}_p := \{ z \in \mathbb{D} : |z| < p \}$. The subclass $C_p$ is given by
 \[
 C_p := \{ f \in \widehat{C_p} : f(0) = 0, \ f^{\prime}(0) = 1 \}.
 \]
  It is important to note that the function 
  \begin{align*}
  	k_p(z) := \frac{pz}{(p-z)(1-pz)}, \quad z \in \mathbb{D}_p,
  \end{align*}
  is a classical example of a concave univalent function with a pole at $p$, whose series expansion is given by 
  $$k_p(z) = \sum_{n=1}^{\infty} \frac{1-p^{2n}}{(1-p^2)p^{n-1}} z^n,$$ and is commonly used as an extremal function.
  \vspace{2mm}
 
 The representation formula for functions in $C_p$ was established by Wirths \cite{Wirths-SMJ-2006} in 2006. This formulation allowed Bhowmik \emph{et al.} \cite{Bhowmik-Ponnusamy-Wirths-KMJ-2007} to initiate an analysis in 2007, based on the Laurent series expansion around the pole $p$, where they obtained certain coefficient bounds for the class $\widehat{C_p}$. The class has since been extensively investigated in \cite{Avkhadiev-Pommerenke-Wirths-MN-2004,Avkhadiev-Wirths-FM-2007,Avkhadiev-Wirths-CVEE-2007,Bhowmik-Das-JMAA-2018,Wirths-SMJ-2003}.
 \vspace{2mm}
 
In the following, the sharp Bohr-type inequality is obtained in terms of convex combinations for harmonic mappings whose analytic part is subordinate to functions belonging to the class $\widehat{C_p}$, where $p \in (0,1)$.

\begin{thm}\label{thm-3.1}
Suppose that $f(z)=h(z)+\overline{g(z)}=\sum_{n=0}^{\infty}a_nz^n+\overline{\sum_{n=1}^{\infty}b_nz^n}$ is a sense-preserving $K$-quasiconformal harmonic mapping in $\mathbb{D}$ and $h\prec\phi$, where $\phi\in\widehat{C_p}$ and $p\in(0,1)$. Then for arbitrary $t\in[0,1]$ and $m\in\mathbb{N}$, the inequality
\begin{align*}
	t|h(z^m)|+(1-t)\sum_{n=0}^{\infty}(|a_n|+|b_n|)r^n\leq|\phi(0)|+ d(\phi(0),\partial\phi(\mathbb{D}))
\end{align*}
for $r\leq \min\{R_{K,m,p,t},1/3\}$,
where $R_{K,m,p,t}\in(0,p)$ is the unique root of the equation 
\begin{align}\label{Eqn-3.1}
	tk_p(r^m)+\left(\frac{2K(1-t)}{K+1}\right)k_p(r)-\frac{p}{(1+p)^2}=0.
\end{align}
The radius $R_{K,m,p,t}$ is sharp if $R_{K,m,p,t}\leq1/3$.
\end{thm}
By choosing $K = 1$ in Theorem \ref{thm-3.1}, we obtain the following corollary.
\begin{cor}\label{Cor-3.1}
	Let $f(z)=\sum_{n=0}^{\infty}a_nz^n$ be a an analytic function defined on $\mathbb{D}$ and $f\prec\phi$, where $\phi\in\widehat{C_p}$ and $p\in(0,1)$. Then for arbitrary $t\in[0,1]$ and $m\in\mathbb{N}$, we have 
	\begin{align*}
		t|f(z^m)|+(1-t)\sum_{n=0}^{\infty}|a_n|r^n\leq |\phi(0)|+d(\phi(0),\partial\phi(\mathbb{D}))
	\end{align*}
	for $r\leq \min\{1/3,R_{m,p,t}\}$, where $R_{m,p,t}$ is the unique root in $(0,1)$ of the equation 
	\begin{align*}
	tk_p(r^m)+(1-t)k_p(r)-\frac{p}{(1+p)^2}=0.
	\end{align*}
	The result is sharp if $R_{m,p,t}\leq1/3$.
\end{cor}
\begin{rem}
\begin{enumerate}
	\item [(i)] By choosing $K=1$ (so $k=0$) and $t=0$ in Theorem \ref{thm-3.1}, we obtain
	\begin{align*}
		R_{0,m,p,0}=(1+1/p+p)-(\sqrt{p}+1/\sqrt{p})\sqrt{p+1/p},
	\end{align*}
	which represents the root of the quadratic $pr^2 - 2(1+p+p^2)r + p = 0$ in $(0,p)$. Therefore, the result of \cite[Corollary 1]{Bhowmik-Das-JMAA-2018} is included as a special case of Theorem \ref{thm-3.1}.\vspace{1.2mm}
	
	\item [(ii)] If we take the limit $p \to 1$ in Corollary \ref{Cor-3.1}, this reflects the case where $\phi$ has a pole on $\partial \mathbb{D}$, meaning that $\phi$ reduces to an analytic univalent function in $\mathbb{D}$. Consequently, Corollary \ref{Cor-3.1} coincides with the result in \cite[Theorem 2.3]{Arora-Vinayak-CVVE-2025} for $n=1$.\vspace{1.2mm}
	
	\item [(iii)] If we let $p \to 1$ and $t=0$ in Corollary \ref{Cor-3.1}, then it reduces to the result of \cite[Theorem 1]{Abu-Muhanna-CVEE-2010}.

\end{enumerate}
	
\end{rem}
\begin{proof}[\bf Proof of Theorem \ref{thm-3.1}]
	For a function $\phi$ in $\widehat{C_p}$ represented as $\phi(z)=\sum_{n=0}^{\infty}c_nz^n$, define the Koebe transform for $a \in \mathbb{D} \setminus \{p\}$ by
	\begin{align*}
		F(z) = \frac{\phi\left( \frac{z+a}{1+\bar{a}z} \right) - \phi(a)}{(1 - |a|^2)\phi^{\prime}(a)}.
	\end{align*}
Since $e^{-it}F(ze^{it})\in C_{\left|\frac{p-a}{1-\bar{a}p}\right|}$ for some $t\in\mathbb{R}$, we apply the result from \cite[Chapter 15, p. 137]{Avkhadiev-Wirths-2009} to obtain
\begin{align}\label{eq-3.1}
	d(\phi(0),\partial\phi(\mathbb{D}))\geq\frac{p}{(1 + p)^2} |\phi^{\prime}(0)|.
\end{align}
	
	Since $F(z)=(\phi(z)-\phi(0))/\phi^{\prime}(0)\in C_p$, we obtain (see \cite[p. 26]{Jenkins-MMJ-1962})
	\begin{align}\label{eq-3.2}
		|c_n|\leq|\phi^{\prime}(0)|\frac{1-p^{2n}}{(1-p^2)p^{n-1}},\;\;n\geq 1.
	\end{align}
	Since $h$ is subordinate to $\phi$ and $\phi \in \widehat{C_p}$ with $p \in (0,1)$, the application of Lemma \ref{lem-2.1} together with inequality \eqref{eq-3.2} yields
	
	\begin{align}\label{eq-3.3}
		\sum_{n=1}^{\infty}|a_n|r^n\leq\sum_{n=1}^{\infty}|c_n|r^n&\leq|\phi^{\prime}(0)|\sum_{n=1}^{\infty} \frac{1-p^{2n}}{(1-p^2)p^{n-1}}r^n\\&=|\phi^{\prime}(0)|k_p(r)\;\;\mbox{for}\;\; |z|=r\leq 1/3\nonumber.
	\end{align}
	
	Given that $f$ is $K$-quasiconformal sense-preserving harmonic mapping on $\mathbb{D}$, Schwarz’s Lemma implies the dilatation $\omega = g'/h'$ is analytic with $|\omega(z)| \leq k$ in $\mathbb{D}$, where $K = (1 + k)/(1 - k)$ and $k \in [0,1)$. Therefore, by Lemma \ref{lem-2.2}, the inequality
	\begin{align}\label{eq-3.4}
		\sum_{n=1}^{\infty} |b_n| r^n \leq k \sum_{n=1}^{\infty} |a_n| r^n \leq k |\phi'(0)| k_p(r)
	\end{align}
	holds for $|z| = r \leq \frac{1}{3}$.
	
	By using the subordination $h \prec \phi$ where $\phi \in \widehat{C_p}$ and $p \in (0,1)$, we derive the following inequality
	\begin{align}\label{eq-3.5}
		|h(z^m)| &\leq |h(0)| + |h(z^m) - h(0)| \\\nonumber
		&\leq |\phi(0)| + |\phi(z^m) - \phi(0)| \\\nonumber
		&= |\phi(0)| + \sum_{n=1}^{\infty} |c_n| r^{mn} \\\nonumber
		&\leq |\phi(0)| + |\phi^{\prime}(0)|k_p(r^m),
	\end{align}
	which holds for $|z| = r \leq \frac{1}{3}$.\vspace{2mm}

	Making use of equations \eqref{eq-3.3} to \eqref{eq-3.5} and the inequality \eqref{eq-3.1}, we find
	\begin{align}\label{eq-3.6}
		A_1(r):=&t|h(z^m)|+(1-t)\sum_{n=0}^{\infty}(|a_n|+|b_n|)r^n\\\nonumber\leq&|\phi(0)|+t|\phi^{\prime}(0)|k_p(r^m)+(1-t)(k+1)|\phi^{\prime}(0)|k_p(r)\\\nonumber\leq&|\phi(0)|+\frac{(1+p)^2}{p}d(\phi(0),\partial\phi(\mathbb{D}))\left(tk_p(r^m)+\left(\frac{2K(1-t)}{K+1}\right)k_p(r)\right)\nonumber.
	\end{align}
	Let
	\begin{align*}
		G^t_{k,m,p}(r):=tk_p(r^m)+\left(\frac{2K(1-t)}{K+1}\right)k_p(r)-\frac{p}{(1+p)^2}.
	\end{align*}
	Then, it can be easily seen that $G^t_{k,m,p}$ is strictly increasing function of $r\in(0,p)$. Also,
	
	\begin{align*}
		G^t_{k,m,p}(0)=-\frac{p}{(1+p)^2}\;\;\mbox{and}\;\;\lim_{r\to p}G^t_{k,m,p}(r)=\infty.
	\end{align*}
	  This implies that the equation $G^t_{K,m,p}(r)=0$ has a unique positive root $R_{K,m,p,t}\in(0,p)$. Equation \eqref{eq-3.6} implies that the inequality $$A_1(r) \leq |\phi(0)| + d(\phi(0), \partial\phi(\mathbb{D}))$$ holds true for all $r \leq\min\{1/3, R_{K,m,p,t}\}$, where $R_{K,m,p,t}$ is the positive root of $G^t_{K,m,p}(r) = 0$.
	  \vspace{2mm}
	
	In order to establish the sharpness for the case $ R_{K,m,p,t}\leq 1/3$, we consider the function $f(z) = h(z) + \overline{g(z)}$ defined on $\mathbb{D}$, where
	\begin{align*}
		h(z) = \phi(z) = k_p(z) \quad \text{and} \quad g(z) = k \lambda k_p(z),
	\end{align*}
	with $|\lambda| = 1$ and $k = (K - 1)/(K + 1)$.
	 It follows from the properties of $k_p$ (see \cite[p. 137]{Avkhadiev-Wirths-2009}) that $$\widehat{C} \setminus k_p(\mathbb{D}) = \left[-\frac{p}{(1 - p)^2}, -\frac{p}{(1 + p)^2}\right].$$ As a result, we see that \begin{align*}
	 	d(\phi(0), \partial \phi(\mathbb{D})) = \frac{p}{(1 + p)^2}.
	 \end{align*} Thus, for $z=r$, we have
	\begin{align*}
		t&|k_p(r^m)|+(1-t)\sum_{n=1}^{\infty}\left(\left|\frac{1-p^{2n}}{(1-p^2)p^{n-1}}\right|+\left|\frac{k\lambda(1-p^{2n})}{(1-p^2)p^{n-1}}\right|\right)r^n\\&=t|k_p(r^m)|+(1-t)(1+k)\sum_{n=1}^{\infty}\frac{1-p^{2n}}{(1-p^2)p^{n-1}}r^n\\&=tk_p(r^m)+\left(\frac{2K(1-t)}{K+1}\right)k_p(r)\\&>\frac{p}{(1+p)^2}\\&=|\phi(0)|+d(\phi(0),\partial\phi(\mathbb{D}))
	\end{align*}
	for $r>R_{K,m,p,t}$, where $R_{K,m,p,t}$ is the positive root of the equation $G^t_{K,m,p}(r)=0$. This implies that $R_{K,m,p,t}$ is the best possible.
\end{proof}
The tabulated values of $R_{K,m,p,t}$ for some choices of $K \geq 1$, $m \in \mathbb{N}$, $p\in(0,1)$, and $t \in [0,1]$ are given in Table~4, and their positions are represented in Figure~4.

\begin{table}[ht]
	\centering
	\begin{tabular}{|l|l|l|l|l|l|l|l|l|l|}
		\hline
		$K$& $m$&$p$& $t$&$R_{K,m,p,t}$ \\
		\hline
		$5$& $3$&$0.9$& $0.2$& $0.1383$ \\
		\hline
		$25$& $15$&$0.7$& $0.6$ & $0.1957$ \\
		\hline
		$50$& $30$&$0.8$&  $0.7$ &$0.2384$\\
		\hline
		$35$& $20$&$0.9 $& $0.77$ &$0.2840$ \\
		\hline
		$40$& $25$& $0.9$ &$0.83$ &$0.3323$\\
		\hline
		
	\end{tabular}
	\caption{$R_{K,m,p,t}$ is the smallest positive root of the equation \eqref{Eqn-3.1} in $(0,p)$}
\end{table}
The following result establishes a Bohr-type inequality involving a parameter for harmonic mappings whose analytic part is subordinate to a function in $\widehat{C_p}$, with $p \in (0,1)$.

\begin{thm}\label{thm-3.2}
	Suppose that $f(z)=h(z)+\overline{g(z)}=\sum_{n=0}^{\infty}a_nz^n+\overline{\sum_{n=1}^{\infty}b_nz^n}$ is a sense-preserving $K$-quasiconformal harmonic mapping in $\mathbb{D}$ and $h\prec\phi$, where $\phi\in\widehat{C_p}$ and $p\in(0,1)$. Then for arbitrary $\lambda>0$ and $m\in\mathbb{N}$,
	\begin{align*}
		|h(z^m)|+\lambda\sum_{n=1}^{\infty}(|a_n|+|b_n|)r^n\leq d(\phi(0),\partial\phi(\mathbb{D})) 
	\end{align*}
	for $r\leq \min\{1/3,R^*_{K,m,p,\lambda}\}$, where $R^*_{K,m,p,\lambda}\in(0,p)$ is the unique root of the equation 
	\begin{align}\label{Eqn-3.8}
		k_p(r^m)+\left(\frac{2K\lambda}{K+1}\right)k_p(r)-\frac{p}{(1+p)^2}=0.
	\end{align}
	The result is sharp if  $R^*_{K,m,p,\lambda}\leq 1/3$.
\end{thm}
Setting $K=1$ in Theorem \ref{thm-3.1} leads to the following corollary.

\begin{cor}\label{Cor-3.2}
	Let $f(z)=\sum_{n=0}^{\infty}a_nz^n$ be a an analytic function defined on $\mathbb{D}$ and $f\prec\phi$, where $\phi\in\widehat{C_p}$ and $p\in(0,1)$. Then for arbitrary $t\in[0,1]$ and $m\in\mathbb{N}$, we have 
	\begin{align*}
		|f(z^m)|+\lambda\sum_{n=1}^{\infty}|a_n|r^n\leq |\phi(0)|+d(\phi(0),\partial\phi(\mathbb{D}))
	\end{align*}
	for $r\leq \min\{1/3,R^*_{m,p,\lambda}\}$, where $R^*_{m,p,\lambda}$ is the unique in $(0,1)$ of the equation 
	\begin{align*}
		k_p(r^m)+\lambda k_p(r)-\frac{p}{(1+p)^2}=0.
	\end{align*}
	The result is sharp if $R^*_{m,p,\lambda}\leq1/3$.
\end{cor}
\begin{rem}
	\begin{enumerate}
		\item [(i)] By applying the conditions $m \to \infty$ and $\lambda=1$ to Corollary \ref{Cor-3.2}, one obtains
		$R^*_{\infty,p,1}=(1+1/p+p)-(\sqrt{p}+1/\sqrt{p})\sqrt{p+1/p}$,
		which coincides with the root of $pr^2 - 2(1+p+p^2)r + p = 0$ lying in $(0,p)$. Thus, Corollary \ref{Cor-3.2} incorporates the result of \cite[Corollary 1]{Bhowmik-Das-JMAA-2018} as a special case.\vspace{1.2mm}
		
		\item [(ii)] Taking the limit $p \to 1$ in Corollary \ref{Cor-3.2} gives the case in which $\phi$ has a pole on $\partial \mathbb{D}$, and thus $\phi$ reduces to an analytic univalent function on $\mathbb{D}$. Therefore, Corollary \ref{Cor-3.2} coincides with \cite[Theorem 2.5]{Arora-Vinayak-CVVE-2025} for $n=1$.\vspace{1.2mm}
		
		\item [(iii)] If we take the limits $p \to 1$, $m \to \infty$, and $\lambda$ in Corollary \ref{Cor-3.2}, it reduces to the result of \cite[Theorem 1]{Abu-Muhanna-CVEE-2010}.\vspace{1.2mm}
		
		\item [(iv)] If we let $p \to 1$ and set $\lambda = m = 1$ in Corollary \ref{Cor-3.2}, we obtain $R^*_{1,1,1}=5-2\sqrt{6}$, which matches the conclusion of \cite[Remark 3]{Kayumov-Khammatova-Ponnusamy-JMAA-2021}.
		
	\end{enumerate}
\end{rem}
\begin{proof}[\bf Proof of Theorem \ref{thm-3.2}]
	By following a line of reasoning similar to that used in the proof of Theorem \ref{thm-3.1}, and by applying Lemmas \ref{lem-2.1} and \ref{lem-2.2} along with equation \eqref{eq-3.2}, we derive the inequalities \eqref{eq-3.3}, \eqref{eq-3.4}, and \eqref{eq-3.5}.
	
	 Therefore, by using inequalities \eqref{eq-3.3}, \eqref{eq-3.4}, and \eqref{eq-3.5}, we have
	\begin{align}\label{eq-3.7}
		A_2(r):=&|h(z^m)|+\lambda\sum_{n=1}^{\infty}(|a_n|+|b_n|)r^n\\\nonumber\leq&|\phi(0)|+|\phi^{\prime}(0)|k_p(r^m)+\lambda(k+1)|\phi^{\prime}(0)|k_p(r)\\\nonumber\leq&|\phi(0)|+\frac{(1+p)^2}{p}d(\phi(0),\partial\phi(\mathbb{D}))\left(k_p(r^m)+\left(\frac{2K\lambda}{K+1}\right)k_p(r)\right)\nonumber.\nonumber
	\end{align}
	Let
	\begin{align*}
		H^{\lambda}_{K,m,p}(r):=k_p(r^m)+\left(\frac{2K\lambda}{K+1}\right)k_p(r)-\frac{p}{(1+p)^2}.
	\end{align*}
	It is clear that $	H^{\lambda}_{K,m,p}(r)$ increases strictly with $r$ over the interval $(0, p)$. Moreover,
	
	\begin{align*}
		H^{\lambda}_{K,m,p}(0) = -\frac{p}{(1+p)^2} \quad \text{and} \quad \lim_{r \to p} H^{\lambda}_{K,m,p}(r) = \infty.
	\end{align*}

	 Thus, the equation $H^{\lambda}_{K,m,p}(r) = 0$ has exactly one positive root $R^*_{K,m,p,\lambda}$ in $(0, p)$. Consequently, in view of \eqref{eq-3.7}, the inequality 
	 \begin{align*}
	 	A_2(r) \leq |\phi(0)| + d(\phi(0), \partial\phi(\mathbb{D}))
	 \end{align*}
	 holds for $r \leq \min\{1/3,R^*_{K,m,p,\lambda}\}$.\vspace{2mm}
	
	To establish the sharpness of the result in the case $R^*_{K,m,p,\lambda}\leq 1/3$, we take the function $f(z)=h(z)+\overline{g(z)}$ in $\mathbb{D}$, where
	\begin{align*}
		h(z)=\phi(z)=k_p(z)\quad\text{and}\quad g(z)=k\lambda k_p(z),
	\end{align*}
	with $|\lambda| = 1$ and $k =(K-1)/(K+1)$. It is known that $d(\phi(0),\partial\phi(\mathbb{D}))=p/(1+p)^2$. Hence, for $z=r$, we obtain
	
	\begin{align*}
		|k_p(r^m)|+\lambda	&\sum_{n=1}^{\infty}\left(\left|\frac{1-p^{2n}}{(1-p^2)p^{n-1}}\right|+\left|\frac{k\lambda(1-p^{2n})}{(1-p^2)p^{n-1}}\right|\right)r^n\\&=|k_p(r^m)|+\lambda(1+k)\sum_{n=1}^{\infty}\frac{1-p^{2n}}{(1-p^2)p^{n-1}}r^n\\&=k_p(r^m)+\left(\frac{2K\lambda}{K+1}\right)k_p(r)\\&>\frac{p}{(1+p)^2}\\&=|\phi(0)|+d(\phi(0),\partial\phi(\mathbb{D}))
	\end{align*}
	for $r>R^*_{K,m,p,\lambda}$, where $R^*_{K,m,p,\lambda}$ is the positive root of the equation $H^{\lambda}_{K,m,p}(r)=0$. This shows that $R^*_{K,m,p,\lambda}$ is the best possible.
\end{proof}
In Table~5, the values of $R_{K,m,p,t}$ are tabulated for specific $K \geq 1$, $m \in \mathbb{N}$, $p \in (0,1)$, and $\lambda > 0$, and their locations are represented in Figure~5.

\begin{table}[ht]
	\centering
	\begin{tabular}{|l|l|l|l|l|l|l|l|l|l|}
		\hline
		$K$& $m$&$p$& $\lambda$&$R^*_{K,m,p,\lambda}$ \\
		\hline
		$25$& $20$&$0.7$& $1$& $0.1003$ \\
		\hline
		$25$& $15$&$0.5$& $0.5$ & $0.1498$ \\
		\hline
		$45$& $25$&$0.8$&  $0.35$ &$0.2171$\\
		\hline
		$35$& $20$&$0.6 $& $0.2$ &$0.2738$ \\
		\hline
		$2.6$& $5$& $0.7$ &$0.2$ &$0.3323$\\
		\hline
	
	\end{tabular}
	\caption{$R^*_{K,m,p,\lambda}$ is the smallest root of the equation \eqref{Eqn-3.8} in $(0,p)$}
\end{table}
\section{\bf{Bohr-type inequalities for the family of concave univalent functions}}

Let us consider the class of concave functions, defined as conformal mappings from the unit disk $\mathbb{D}$ onto the complement of convex sets with opening angle $\pi\alpha$ at infinity, for $\alpha \in [1, 2]$. The name ``concave'' reflects the structure of the image domain $f(\mathbb{D})$. It is crucial to note that the set $\mathbb{C} \cup \{\infty\} \setminus f(\mathbb{D})$ can be either bounded or unbounded. If it is bounded, then $f$ is a meromorphic function. Therefore, we restrict ourselves to the unbounded case, where $f$ has a simple pole on the boundary of $\mathbb{D}$. We proceed to define this class formally.

Let $\widehat{\mathcal{C}_0}(\alpha)$ represent the class of analytic functions in $\mathbb{D}$ that fulfill the following conditions:
\begin{enumerate}
	\item $f$ maps $\mathbb{D}$ onto $f(\mathbb{D})$ conformally.
	\item $f(1) = \infty$.
	\item The set $\mathbb{C} \cup \{\infty\} \setminus f(\mathbb{D})$ is convex, and the opening angle at infinity is less than or equal to $\pi\alpha$, where $\alpha \in [1, 2]$.
\end{enumerate}
This collection of functions is known as the class of concave univalent functions with opening angle $\pi\alpha$. The function $f_\alpha$, defined on $\mathbb{D}$ as
\begin{align}\label{eq-4.1}
	f_\alpha(z):=\frac{1}{2\alpha}\left[\left(\frac{1+z}{1-z} \right)^\alpha-1\right]
	= \sum_{k=1}^\infty A_kz^k.
\end{align}
is of great significance in this class and commonly acts as an extremal function. It is easy to verify that for all $\alpha \in [1, 2]$, the class $\widehat{\mathcal{C}_0}(\alpha)$ is contained in $\widehat{\mathcal{C}_0}(2)$ (see  \cite{Avkhadiev-Wirths-LJM-2005}). In the special case $\alpha = 1$, the image domain is a convex half-plane. This shows that concave univalent functions relate to convex functions, and every $f \in\widehat{\mathcal{C}_0}(1)$ belongs to the class of convex functions. For a more comprehensive analysis of concave univalent functions, we refer the reader to \cite{Avkhadiev-Pommerenke-Wirths-CMH-2006, Avkhadiev-Wirths-LJM-2005,Bhowmik-MN-2012,Bhowmik-Das-JMAA-2018,Cruz-Pommerenke-CVEE-2007} and the works cited therein.\vspace{2mm}

In order to prove the main results, we rely on the following lemma from \cite{Bhowmik-MN-2012}. This lemma provides both the Growth \cite[Corollary 2.4]{Bhowmik-MN-2012} and Distortion \cite[Corollary 2.3]{Bhowmik-MN-2012} theorems for concave univalent functions.
\begin{lem}
	Let $\alpha\in[1,2]$. Then, for each $f\in\widehat{C_0}(\alpha)$, we have 
	\begin{enumerate}
		\item [(i)] Growth Theorem: $|f(z)-f(0)|\leq |f^{\prime}(0)|f_{\alpha}(r)$,\vspace{1.2mm}
		
		\item[(ii)] Distortion Theorem:  $|f^{\prime}(z)|\leq |f^{\prime}(0)|f^{\prime}_{\alpha}(r)$,
	\end{enumerate}
	where $f_{\alpha}$ is defined by \eqref{eq-4.1}. For each $z\in\mathbb{D}$, $z\neq0$, equality occurs if $f=f_{\alpha}$.
\end{lem}
\begin{lem}\cite{Allu-Arora-JMAA-2023}
	Let $f\in\widehat{C_0}(\alpha)$, $\alpha\in[1,2]$, have the expansion $f(z)=\sum_{n=0}^{\infty}a_nz^n$. Then we have the following inequalities:
	\begin{enumerate}
		\item [(i)] $|f^{\prime}(0)|\leq 2\alpha d(f(0), \partial f(\mathbb{D}))$.\vspace{1.2mm}
		
		\item[(ii)] $|a_n|\leq A_n|f^{\prime}(0)|$, for $n\geq 1$. 
	\end{enumerate}
	All inequalities are sharp for the function $f_{\alpha}(z)$.
\end{lem}
The following theorem provides a Bohr-type inequality in terms of convex combinations for harmonic mappings with analytic part subordinate to functions in the class $\widehat{C_0}(\alpha)$, where $\alpha \in [1,2]$.

\begin{thm}\label{thm-4.1}
	Suppose that $f(z)=h(z)+\overline{g(z)}=\sum_{n=0}^{\infty}a_nz^n+\overline{\sum_{n=1}^{\infty}b_nz^n}$ is a sense-preserving $K$-quasiconformal harmonic mapping in $\mathbb{D}$ and $h\prec\phi$, where $\phi\in\widehat{C_0}(\alpha)$ and $\phi(z)=\sum_{n=0}^{\infty}c_nz^n$. Then for arbitrary $t\in[0,1]$ and $m\in\mathbb{N}$, the inequality
	\begin{align*}
		t|h(z^m)|+(1-t)\sum_{n=0}^{\infty}(|a_n|+|b_n|)r^n\leq|\phi(0)|+ d(\phi(0),\partial\phi(\mathbb{D})),
	\end{align*}
	holds for $|z|=r\leq\min\{1/3, R_{K,m,t,\alpha}\}$, where $R_{K,m,t,\alpha}\in(0,1)$ is the unique root of the equation 
	\begin{align}\label{Eqn-4.2}
		tf_{\alpha}(r^m)+\left(\frac{2K(1-t)}{K+1}\right)f_{\alpha}(r)-\frac{1}{2\alpha}=0.
	\end{align}
	The radius $R_{K,m,t,\alpha}$ is sharp if $R_{K,m,t,\alpha}\leq 1/3$.
\end{thm}
\begin{rem}
	\begin{enumerate}
		\item [(i)] By choosing $K=1$ in Theorem \ref{thm-4.1}, we obtain the result stated in \cite[Theorem 4.2]{Arora-Vinayak-CVVE-2025}.\vspace{1.2mm}
		
		\item[(ii)]  Choosing $K=1$ and $t=0$ in Theorem \ref{thm-4.1} leads to the result established in \cite{Bhowmik-Das-JMAA-2018}.
	\end{enumerate}
\end{rem}
\begin{proof}[\bf Proof of Theorem \ref{thm-4.1}]
	Since $h\prec\phi$ and $\phi\in\widehat{C_0}(\alpha)$, by applying Lemmas 4.1(a) and 4.1(d), we have
	\begin{align}\label{Eq-3.1}
	|h(z^m)|\leq&|h(0)|+|h(z^m)-h(0)|\\\nonumber\leq&|\phi(0)|+|\phi(z^m)-\phi(0)|\\\nonumber\leq& |\phi(0)|+2\alpha d(\phi(0),\partial \phi(\mathbb{D}))f_{\alpha}(r^m).
	\end{align}
	Using Lemma 2.2(a) and the inequality from Lemma 4.1(c), it follows that 
	\begin{align*}
		\sum_{n=0}^{\infty}|a_n|r^n\leq \sum_{n=0}^{\infty}|c_n|r^n\leq |\phi(0)|+\phi^{\prime}(0)|f_{\alpha}(r)\;\;\mbox{for}\;\; r\leq1/3.
	\end{align*}
	With another application of Lemma 4.1(d) to the inequality above, we find 
	\begin{align}\label{Eq-3.2}
		\sum_{n=0}^{\infty}|a_n|r^n\leq |\phi(0)|+2\alpha d(\phi(0),\partial \phi(\mathbb{D}))f_{\alpha}(r),\;\; r\leq\frac{1}{3}.
	\end{align}

By applying a similar reasoning as in the proof of Theorem \ref{thm-2.1} and considering Lemma \ref{lem-2.2}, we obtain
	\begin{align}\label{Eq-3.3}
	\sum_{n=1}^{\infty}|b_n|r^n\leq k\sum_{n=1}^{\infty}|a_n|r^n\leq 2k\alpha d(\phi(0),\partial \phi(\mathbb{D}))f_{\alpha}(r),\;\; r\leq\frac{1}{3}.
\end{align}
Hence, from the combined use of \eqref{Eq-3.1}, \eqref{Eq-3.2}, and \eqref{Eq-3.3}, we deduce that
\begin{align}\label{eq-4.5}
B_1(r):=&t|h(z^m)|+(1-t)\sum_{n=0}^{\infty}(|a_n|+|b_n|)r^n\\\nonumber\leq& |\phi(0)|+2\alpha d(\phi(0),\partial \phi(\mathbb{D}))\left(tf_{\alpha}(r^m)+\left(\frac{2K(1-t)}{K+1}\right)f_{\alpha}(r)\right)\\\nonumber=&|\phi(0)|+2\alpha d(\phi(0),\partial \phi(\mathbb{D}))\left(F_{k,m,t,\alpha}(r)+\frac{1}{2\alpha}\right), 
\end{align}
where $$F_{K,m,t,\alpha}(r)=tf_{\alpha}(r^m)+\left(\frac{2K(1-t)}{K+1}\right)f_{\alpha}(r)-\frac{1}{2\alpha}.$$ It is straightforward to verify that $F_{K,m,t,\alpha}(r)$ is continuous and strictly increasing on the interval $(0,1)$. We observe that
\begin{align*}
	F_{K,m,t,\alpha}(0)<0\;\;\;\mbox{and}\;\;\lim_{r\to1}F_{k,m,t,\alpha}(r)=+\infty,
\end{align*}
 which guarantees a unique root $R_{K,m,t,\alpha}$ of the equation $F_{K,m,t,\alpha}(r)=0$ in $(0,1)$. As a result, the inequality $F_{K,m,t,\alpha}(r)\leq 0$ holds for all $r \leq R_{K,m,t,\alpha}$, which in turn implies from \eqref{eq-3.4} that 
 \begin{align*}
 	B_1(r)\leq |\phi(0)|+d(\phi(0),\partial\phi(\mathbb{D}))\;\;\mbox{for}\;\;r \leq \min\{1/3, R_{K,m,t,\alpha}\}
 \end{align*}\vspace{2mm}

To prove that the result is sharp for the case $R_{k,m,t,\alpha}\leq 1/3$, we consider the function $f(z) = h(z) + \overline{g(z)}$ in the unit disk $\mathbb{D}$, where
\begin{align*}
	h(z) = \phi(z) = f_{\alpha}(z) \quad \text{and} \quad g(z) = k\lambda f_{\alpha}(z)
\end{align*}
with $|\lambda|=1$ and $k =(K-1)/(K+1)$. It is known that $d(f_{\alpha}(0), f_{\alpha}(\mathbb{D})) = 1/2\alpha$. Hence, for $z=r$, we have
\begin{align*}
	t|f_{\alpha}(r^m)|+(1-t)\sum_{n=1}^{\infty}(|A_n|+|\lambda k A_n|)r^n=&t|f_{\alpha}(r^m)|+\left(\frac{2K(1-t)}{K+1}\right)f_{\alpha}(r)\\>\frac{1}{2\alpha}=|\phi(0)|+d(\phi(0),\partial\phi(\mathbb{D}))
\end{align*}
for $r>R_{K,m,t,\alpha}$. This completes the proof.

\end{proof}
The tabulated values of $R_{K,m,t,\alpha}$ for $K \geq 1$, $m \in \mathbb{N}$, $t \in [0,1]$, and $\alpha \in [1,2]$ are given in Table~6, along with their visualization in Figure~6.  
\begin{table}[ht]
	\centering
	\begin{tabular}{|l|l|l|l|l|l|l|l|l|l|}
		\hline
		$K$& $m$&$t$& $\alpha$&$R_{K,m,t,\alpha}$ \\
		\hline
		$50$& $7$&$0.2$& $2$& $0.1227$ \\
		\hline
		$30$& $10$&$0.3$& $1.5$ & $0.1822$ \\
		\hline
		$7$& $15$&$0.5$&  $1.6$ &$0.2338$\\
		\hline
		$45$& $9$&$0.5 $& $1.2$ &$0.2853$ \\
		\hline
		$1$& $1$& $1$ &$1$ &$1/3$\\
		\hline
	\end{tabular}
	\caption{$R_{K,m,t,\alpha}$ is the unique positive root of the equation \eqref{Eqn-4.2} in $(0,1)$}
\end{table}
The following theorem establishes a Bohr-type inequality with a parameter for harmonic mappings whose analytic part is subordinate to functions belonging to the class $\widehat{C_0}(\alpha)$, where $\alpha \in [1,2]$.

\begin{thm}\label{thm-4.2}
		Suppose that $f(z)=h(z)+\overline{g(z)}=\sum_{n=0}^{\infty}a_nz^n+\overline{\sum_{n=1}^{\infty}b_nz^n}$ is a sense-preserving $K$-quasiconformal harmonic mapping in $\mathbb{D}$ and $h\prec\phi$, where $\phi\in\widehat{C_0}(\alpha)$ and $\phi(z)=\sum_{n=0}^{\infty}c_nz^n$. Then for arbitrary $\lambda>0$ and $m\in\mathbb{N}$, the inequality
	\begin{align*}
		|h(z^m)|+\lambda\sum_{n=1}^{\infty}(|a_n|+|b_n|)r^n\leq|\phi(0)|+ d(\phi(0),\partial\phi(\mathbb{D})),
	\end{align*}
	holds for $|z|=r\leq\min\{1/3, R^*_{K,m,\alpha,\lambda}\}$, where $R^*_{K,m,\alpha,\lambda}\in(0,1)$ is the unique root of the equation 
	\begin{align}\label{Eqn-4.7}
		f_{\alpha}(r^m)+\left(\frac{2K\lambda}{K+1}\right)f_{\alpha}(r)-\frac{1}{2\alpha}=0.
	\end{align}
	The radius $R^*_{K,m,\alpha,\lambda}$ is sharp if $R^*_{K,m,\alpha,\lambda}\leq 1/3$.
\end{thm}
\begin{rem}
If we set $K=1$ in Theorem \ref{thm-4.2}, the result coincides with \cite[Theorem 4.4]{Arora-Vinayak-CVVE-2025}.  
\end{rem}
\begin{proof}[\bf Proof of Theorem \ref{thm-4.2}]
Following the same line of reasoning as in the proof of Theorem \ref{thm-3.1}, we find that for $r \leq 1/3$, inequalities \eqref{Eq-3.1}, \eqref{Eq-3.2}, and \eqref{Eq-3.3} imply
\begin{align}\label{eq-4.6}
	B_2(r):=&|h(z^m)|+\lambda\sum_{n=0}^{\infty}(|a_n|+|b_n|)r^n\\\nonumber\leq& |\phi(0)|+2\alpha d(\phi(0),\partial \phi(\mathbb{D}))\left(f_{\alpha}(r^m)+\left(\frac{2K\lambda}{K+1}\right)f_{\alpha}(r)\right)\\\nonumber=&|\phi(0)|+2\alpha d(\phi(0),\partial \phi(\mathbb{D}))\left(G_{k,m,\alpha,\lambda}(r)+\frac{1}{2\alpha}\right), 
\end{align}
where
\begin{align*}
	G_{K,m,\alpha,\lambda}(r)=f_{\alpha}(r^m)+\left(\frac{2K\lambda}{K+1}\right)f_{\alpha}(r)-\frac{1}{2\alpha}.
\end{align*}
 It is easy to confirm that $G_{K,m,\alpha,\lambda}(r)$ is continuous and strictly increasing on $(0,1)$. Since 
\begin{align*}
	G_{K,m,\alpha,\lambda}(0)<0\;\;\mbox{and}\;\;\lim_{r\to 1} G_{K,m,\alpha,\lambda}(r) = +\infty,
\end{align*}
 the equation $G_{K,m,\alpha,\lambda}(r)=0$ admits a unique root $R^*_{K,m,\alpha,\lambda}$ in $(0,1)$. Consequently, $G_{K,m,\alpha,\lambda}(r)\leq 0$ for all $r \leq R^*_{K,m,\alpha,\lambda}$, and it follows from \eqref{eq-3.5} that \begin{align*}
 	B_2(r)\leq |\phi(0)| + d(\phi(0), \partial\phi(\mathbb{D}))\;\;\mbox{for}\;\;r \leq \min\{1/3, R^*_{K,m,\alpha,\lambda}\}.
 \end{align*}

To verify the sharpness of the result in the case $R^*_{K,m,\alpha,\lambda}\leq 1/3$, we take the function $f(z) = h(z) + \overline{g(z)}$ defined on the unit disk $\mathbb{D}$, where
\[
h(z) = \phi(z) = f_{\alpha}(z) \quad \text{and} \quad g(z) = k\mu f_{\alpha}(z),
\]
with $|\mu|=1$ and $k =(K-1)/(K+1)$. It is well known that $d(f_{\alpha}(0), f_{\alpha}(\mathbb{D})) =1/2\alpha$. Therefore, for $z = r$, we obtain
\begin{align*}
	|f_{\alpha}(r^m)|+\lambda\sum_{n=1}^{\infty}(|A_n|+|\mu k A_n|)r^n
	= &f_{\alpha}(r^m)+\left(\frac{2K\lambda}{K+1}\right)f_{\alpha}(r)\\
	>& \frac{1}{2\alpha} = |\phi(0)|+d(\phi(0),\partial\phi(\mathbb{D}))
\end{align*}
for $r > R^*_{K,m,\alpha, \lambda}$. This completes the proof.
\end{proof}
The values of $R^*_{K,m,\alpha,\lambda}$ are displayed in Table~7 for $K \geq 1$, $m \in \mathbb{N}$, $\alpha \in [1,2]$, $\lambda>0$, and their visualization is provided in Figure~7.
\begin{table}[ht]
	\centering
	\begin{tabular}{|l|l|l|l|l|l|l|l|l|l|}
		\hline
		$K$& $m$& $\alpha$& $\lambda$&$R^*_{K,m,\alpha,\lambda}$ \\
		\hline
		$5$& $10$&$1.1$& $2$& $0.1187$ \\
		\hline
		$40$& $20$&$2$& $0.5$ & $0.1746$ \\
		\hline
		$30$& $15$&$1.2$&  $0.7$ &$0.2263$\\
		\hline
		$20$& $12$&$1.5 $& $0.4$ &$0.2724$ \\
		\hline
		$25$& $13$& $1.5$ &$0.3$ &$0.3232$\\
		\hline
		
	\end{tabular}
	\caption{$R^*_{K,m,\alpha,\lambda}$ is the unique positive root of the equation \eqref{Eqn-4.7} in $(0,1)$}
\end{table}
In the following result, by replacing $a_1$ with $h^{\prime}(z^m)$, we establish a similar inequality involving a parameter for harmonic mappings whose analytic part is subordinate to a function in $\widehat{C_0}(\alpha)$, where $\alpha \in [1,2]$.

\begin{thm}\label{thm-4.3}
	Suppose that $f(z)=h(z)+\overline{g(z)}=\sum_{n=0}^{\infty}a_nz^n+\overline{\sum_{n=1}^{\infty}b_nz^n}$ is a sense-preserving $K$-quasiconformal harmonic mapping in $\mathbb{D}$ and $h\prec\phi$, where $\phi\in\widehat{C_0}(\alpha)$ and $\phi(z)=\sum_{n=0}^{\infty}c_nz^n$. Then for arbitrary $\lambda>0$ and $m\in\mathbb{N}$, we have 
	\begin{align*}
		|h(z^m)|+|h^{\prime}(z^m)|r+\lambda\left(\sum_{n=2}^{\infty}|a_n|r^n+\sum_{n=1}^{\infty}|b_n|r^n\right)\leq|\phi(0)|+ d(\phi(0),\partial\phi(\mathbb{D})),
	\end{align*}
 for $|z|=r\leq\min\{1/3, R^{**}_{K,m,\alpha,\lambda}\}$, where $R^{**}_{K,m,\alpha,\lambda}\in(0,1)$ is the unique root of the equation 
	\begin{align*}
		f_{\alpha}(r^m)+r\frac{(1+r^m)^{\alpha-1}}{(1-r^m)^{\alpha+1}}+\lambda\left(\frac{2K}{K+1}f_{\alpha}(r)-r\right)-\frac{1}{2\alpha}=0.
	\end{align*}
	The result is sharp if  $R^{**}_{K,m,\alpha,\lambda}\leq 1/3$.
\end{thm}
\begin{rem}
	Choosing $K=1$ in Theorem \ref{thm-4.3} leads directly to the result presented in \cite[Theorem 4.5]{Arora-Vinayak-CVVE-2025}.
\end{rem}
\begin{proof}[\bf Proof of Theorem \ref{thm-4.3}]
	As $h$ is subordinate to $\phi$, there exists a Schwarz function $w$ such that $h(z) = \phi(w(z))$ for all $z \in \mathbb{D}$. Applying the Schwarz–Pick inequality along with Lemma 4.1(b), we obtain
	\begin{align*}
		|h^{\prime}(z^m)|=&|\phi^{\prime}(w(z^m))w^{\prime}(z^m)|\\\leq&|\phi^{\prime}(0)||f^{\prime}_{\alpha}(w(z^m))|\left(\frac{1-|w(z^m)|^2}{1-|z^m|^2}\right)\\=&|\phi^{\prime}(0)|\left(\frac{(1+|w(z^m)|)^{\alpha-1}}{(1-|w(z^m)|)^{\alpha+1}}\right)\left(\frac{1-|w(z^m)|^2}{1-|z^m|^2}\right)\\\leq&|\phi^{\prime}(0)|\frac{(1+r^m)^{\alpha-1}}{(1-r^m)^{\alpha+1}}.
	\end{align*}
	
	Applying Lemma 4.1(d) to the preceding inequality, we obtain
	\begin{align}\label{eq-4.7}
		|h^{\prime}(z^m)|\leq 2\alpha d(\phi(0),\partial\phi(\mathbb{D}))\frac{(1+r^m)^{\alpha-1}}{(1-r^m)^{\alpha+1}}.
	\end{align}
	Applying Lemma \ref{lem-2.1} for $N=2$, along with the inequalities from Lemma 4.1(c) and 4.1(d), we obtain
	\begin{align}\label{eq-4.8}
		\sum_{n=2}^{\infty}|a_n|r^n&\leq \sum_{n=2}^{\infty}|c_n|r^n\\\nonumber&\leq|\phi^{\prime}(0)|\sum_{n=2}^{\infty}A_nr^n\\\nonumber&\leq2\alpha d(\phi(0),\partial\phi(\mathbb{D}))(f_{\alpha}(r)-r) \;\;\mbox{for}\;\; r\leq1/3.
	\end{align}
	Using the inequalities \eqref{Eq-3.1}, \eqref{Eq-3.3}, \eqref{eq-4.7}, and \eqref{eq-3.7} for $|z| = r \leq \frac{1}{3}$, we obtain
	\begin{align}
		B_3(r):=&|h(z^m)|+|h^{\prime}(z^m)|r+\lambda\left(\sum_{n=2}^{\infty}|a_n|r^n+\sum_{n=1}^{\infty}|b_n|r^n\right)\\\nonumber\leq&|\phi(0)|+2\alpha d(\phi(0),\partial\phi(\mathbb{D}))\left(f_{\alpha}(r^m)+r\frac{(1+r^m)^{\alpha-1}}{(1-r^m)^{\alpha+1}}+\lambda(f_{\alpha}(r)-r+kf_{\alpha}(r))\right)\\\nonumber=&|\phi(0)|+2\alpha d(\phi(0),\partial\phi(\mathbb{D}))\left(f_{\alpha}(r^m)+r\frac{(1+r^m)^{\alpha-1}}{(1-r^m)^{\alpha+1}}+\lambda\left(\frac{2K}{K+1}f_{\alpha}(r)-r\right)\right)\\\nonumber=&|\phi(0)|+2\alpha d(\phi(0),\partial\phi(\mathbb{D}))\left(H_{K,m,\alpha,\lambda}(r)+\frac{1}{2\alpha}\right),
	\end{align}
	where
	\begin{align*}
	H_{K,m,\alpha,\lambda}(r)=f_{\alpha}(r^m)+r\frac{(1+r^m)^{\alpha-1}}{(1-r^m)^{\alpha+1}}+\lambda\left(\frac{2K}{K+1}f_{\alpha}(r)-r\right)-\frac{1}{2\alpha}.
	\end{align*}
	The function $H_{K,m,\alpha,\lambda}(r)$ is easily seen to be continuous and strictly increasing on $(0,1)$, with
	\begin{align*}
	H_{K,m,\alpha,\lambda}(0)<0\;\;\mbox{and}\;\;\lim_{r\to 1} H_{K,m,\alpha,\lambda}(r)=+\infty,
	\end{align*}
	so the function admits a unique zero $R^{**}_{K,m,\alpha,\lambda}$ in $(0,1)$. This implies $H_{K,m,\alpha,\lambda}(r)\leq 0$ for $r\leq R^{**}_{K,m,\alpha,\lambda}$, and hence by \eqref{eq-3.5}, we have
	\begin{align*}
	B_3(r)\leq |\phi(0)| + d(\phi(0), \partial\phi(\mathbb{D}))\;\;\mbox{for}\;\;r \leq \min\{1/3, R^{**}_{K,m,\alpha,\lambda}\}.
	\end{align*}
	
	To establish the sharpness of the result for the case $R^{**}_{K,m,\alpha,\lambda}\leq 1/3$, we take the function $f(z)=h(z)+\overline{g(z)}$ on the unit disk $\mathbb{D}$, where  
	\begin{align*}
		h(z)=\phi(z)=f_{\alpha}(z) \quad \text{and} \quad g(z)=k\mu f_{\alpha}(z),
	\end{align*}  
	with $|\mu|=1$ and $k=(K-1)/(K+1)$. It is well known that $d(f_{\alpha}(0),f_{\alpha}(\mathbb{D}))=1/2\alpha$. Moreover, by definition of $f_{\alpha}$,  
	\begin{align*}
		f^{\prime}_{\alpha}(z)=\frac{(1+z)^{\alpha-1}}{(1-z)^{\alpha+1}}.
	\end{align*}  
	Hence, for $z = r$, we obtain
	\begin{align*}
		|f_{\alpha}(r^m)|+&|f^{\prime}_{\alpha}(r^m)|r+\lambda\left(\sum_{n=2}^{\infty}|A_n|r^n+\sum_{n=1}^{\infty}|\mu k A_n|r^n\right)\\
		=&f_{\alpha}(r^m)+r\frac{(1+r^m)^{\alpha-1}}{(1-r^m)^{\alpha+1}}+\lambda\left(\frac{2K}{K+1}f_{\alpha}(r)-r\right)\\
		>& \frac{1}{2\alpha}\\=& |\phi(0)|+d(\phi(0),\partial\phi(\mathbb{D}))
	\end{align*}
	for $r > R^*_{K,m,\alpha, \lambda}$. This completes the proof.
\end{proof}

\noindent{\bf Acknowledgment:} The first author is supported by Science and Engineering Research Board (SERB) (File No. SUR/2022/002244), Govt. of India, and the second author is supported by UGC-JRF (NTA Ref. No.: $ 221610103011$), New Delhi, India. \vspace{2mm}

%\section{Declaration}
\noindent\textbf{Compliance of Ethical Standards}\\

\noindent\textbf{Conflict of interest.} The authors declare that there is no conflict  of interest regarding the publication of this paper.\vspace{1.5mm}

\noindent\textbf{Data availability statement.}  Data sharing not applicable to this article as no datasets were generated or analyzed during the current study.\vspace{1.5mm}

\noindent\textbf{Funding.} No fund.


\begin{thebibliography}{200}
	
\bibitem{Abu-Muhanna-CVEE-2010} {\sc Y. Abu-Muhanna}, Bohr’s phenomenon in subordination and bounded harmonic classes, {\it Complex Var. Elliptic Equ.} {\bf 55} (2010), 1071–1078.

\bibitem{Abu-Muhanna-Ali-JMAA-2011} {\sc Y. Abu-Muhanna} and  {\sc R. M. Ali}, Bohr’s phenomenon for analytic functions into the exterior of a compact convex body, {\it J. Math. Anal. Appl.} {\bf 379} (2011), 512–517.

\bibitem{Abu-Muhanna-Ali-Ponnusami-2016} {\sc Y. Abu-Muhanna}, {\sc R. M. Ali}, and {\sc S. Ponnusamy}, On the Bohr inequality, in: N.K. Govil, et al. (Eds.), {\it  Progress in Approximation Theory and Applicable Complex Analysis}, in: Springer Optimization and Its Applications, vol. 117, 2016, pp. 265–295.

\bibitem{Ahamed-CMFT-2022} {\sc M. B. Ahamed}, The Bohr–Rogosinski radius for a certain class of close-to-convex harmonic mappings, {\it Comput. Methods Funct. Theory} (2022), 1-19.

\bibitem{Ahamed-Ahammed-MJM-2024} {\sc M. B. Ahamed} and {\sc S. Ahammed}, Bohr Inequalities for Certain Classes of Harmonic Mappings, {\it Mediterr. J. Math.} {\bf 21}(1) (2024), 21.

\bibitem{Ahamed-Allu-Halder-AFM-2022} {\sc M. B. Ahamed}, {\sc V. Allu}, and {\sc H. Halder}, The Bohr phenomenon for analytic functions on a
shifted disk, {\it Ann. Fenn. Math.} {\bf 47} (2022), 103–120.

\bibitem{Ali-Jain-Ravichandran-RM-2019} {\sc R. M. Ali}, {\sc N. K. Jain}, and {\sc V. Ravichandran}, Bohr radius for classes of analytic
functions, {\it Results Math.} {\bf 74}(4) (2019), 179.

\bibitem{Alkhaleefah-Kayumov-Ponnusamy-PAMS-2019} {\sc S. A. Alkhaleefah}, {\sc I. R. Kayumov}, and {\sc S. Ponnusamy}, On the Bohr inequality with a fixed zero coefficient, {\it Proc. Amer. Math. Soc.} {\bf 147} (2019), 5263–5274.

%\bibitem{Allu-Ghosh-PIAS-2023} {\sc V. Allu} and {\sc N. Ghosh}, Bohr type inequality for Ces$\acute{a}$ro and Bernardi integral operator on simply connected domain, {\it Proc. Indian Acad. Sci. (Math. Sci.)} {\bf 133} (22) (2023).

\bibitem{Allu-Halder-JMAA-2021} {\sc V. Allu} and {\sc H. Halder}, Bohr radius for certain classes of starlike and convex univalent functions,
{\it J. Math. Anal. Appl.} {\bf 493}(1) (2021), 124519.

\bibitem{Allu-Arora-JMAA-2023} {\sc V. Allu} and {\sc V. Arora}, Bohr-Rogosinski type inequalities for concave univalent functions, {\it J. Math. Anal. Appl.} {\bf 520} (1) (2023), 126845.

\bibitem{Anand-Jain-Kumar-BMMSS-2021} {\sc S. Anand}, {\sc N. K. Jain}, and {\sc S. Kumar}, Sharp Bohr radius constants for certain analytic functions, {\it Bull. Malays. Math. Sci. Soc.} {\bf 44} (2021), 1771-1785.

\bibitem{Arora-Vinayak-CVVE-2025} {\sc V. Arora} and {\sc M. Vinayak}, Bohr's phenomenon for certain classes of analytic functions, {\it Complex Var. Elliptic Equ.}, 1–23. https://doi.org/10.1080/17476933.2024.2439960

\bibitem{Avkhadiev-Pommerenke-Wirths-MN-2004} {\sc F. G. Avkhadiev}, {\sc Ch. Pommerenke}, and {\sc K. -J. Wirths}, On the coefficients of concave univalent functions, {\it Math. Nachr.} {\bf 271} (2004), 3–9.

\bibitem{Avkhadiev-Pommerenke-Wirths-CMH-2006} {\sc F. G. Avkhadiev}, {\sc Ch. Pommerenke}, and {\sc K. -J. Wirths}, Sharp inequalities for the coefficient of concave schlicht functions, {\it  Com-ment. Math. Helv.} {\bf 81} (2006), 801-807.

\bibitem{Avkhadiev-Wirths-LJM-2005} {\sc F. G. Avkhadiev} and {\sc  K. -J. Wirths}, Concave schlicht functions with bounded opening angle at infinity, {\it Lobachevskii J. Math.} {\bf 17} (2005), 3–10.

\bibitem{Avkhadiev-Wirths-FM-2007} {\sc F. G. Avkhadiev} and {\sc  K. -J. Wirths}, A proof of Livingston conjecture, {\it Forum Math.} {\bf 19} (2007), 149–158.

\bibitem{Avkhadiev-Wirths-CVEE-2007} {\sc F. G. Avkhadiev} and {\sc  K. -J. Wirths}, Subordination under concave univalent functions, {\it Complex Var. Elliptic Equ.} {\bf 52}(4) (2007), 299–305.

\bibitem{Avkhadiev-Wirths-2009} {\sc F. G. Avkhadiev} and {\sc K. -J. Wirths}, Schwarz-Pick Type Inequalities, Frontiers in Mathematics, Birkh$\ddot{a}$user Verlag, Basel, 2009, viii+156 pp.

\bibitem{Bhowmik-MN-2012} {\sc B. Bhowmik}, On concave univalent function, {\it  Math. Nachr.} {\bf 285}(5-6) (2012), 606-612.

\bibitem{Bhowmik-Das-JMAA-2018} {\sc B. Bhowmik} and {\sc N. Das}, Bohr phenomenon for subordinating families of certain univalent
functions, {\it J. Math. Anal. Appl.} {\bf  462}(2) (2018), 1087–1098.

\bibitem{Bhowmik-Ponnusamy-Wirths-KMJ-2007} {\sc B. Bhowmik}, {\sc S. Ponnusamy}, and {\sc K. -J. Wirths}, Domains of variability of Laurent coefficients and the convex hull for the family of concave univalent functions, {\it Kodai Math. J.} {\bf 30} (2007), 385–393.

\bibitem{Bohr-PLMS-1914} {\sc H. Bohr}, A theorem concerning power series, {\it Proc. Lond. Math. Soc.} {\bf s2-13} (1914), 1–5.

\bibitem{Cruz-Pommerenke-CVEE-2007} {\sc L. Cruz} and {\sc Ch. Pommerenke}, On concave univalent functions, {\it Complex Var. Elliptic Equ.} {\bf 52}(2–3) (2007), 153–159.
	
%\bibitem{Brawn-Halmos-ASM-1965}	{\sc A. Brawn }, {\sc P. R. Halmos}, and {\sc A. L. Shields}, Ces$\acute{a}$ro operators, {\it Acta Sci. Math.} (Szeged) {\bf 26} (1965), 125–137.

\bibitem{Djakov-JA-2000} {\sc P. B. Djakov} and {\sc M. S. Ramanujan}, A remark on Bohr’s theorem and its generalizations, \textit{J. Anal.} \textbf{8}(2000), 65--77.


\bibitem{Duren-1983} {\sc P. L. Duren}, Univalent functions, Springer, New York, 1983. 

\bibitem{Duren-2004} {\sc P. L. Duren}, Harmonic mapping in the plane, Cambridge University Press (2004).

 \bibitem{Evdoridis-Ponnusamy-Rasila} {\sc S. Evdoridis}, {\sc S. Ponnusamy}, and {\sc A. Rasila}, Improved Bohr’s inequality for shifted disks, {\it  Results Math.} {\bf 76} (14) (2021).
	
\bibitem{Fournier-Ruscheweyh-CRM-2010} {\sc R. Fournier} and {\sc St. Ruscheweyh}, On the Bohr radius for a simply connected plane domain,
Centre de Recherches Mathematiques CRM Proceeding and Lecture Notes, {\bf 51} (2010), 165–171.

\bibitem{Gangania-Kumar-MJM-2022} {\sc K. Gangania} and {\sc S. S. Kumar}, Bohr-Rogosinski phenomenon for $\mathcal{S}^*(\psi)$ and $\mathcal{C}(\psi)$, {\it Mediterr. J. Math.} {\bf 19} (2022), 161.

\bibitem{Garcia-Mashreghi-Ross-2018} {\sc S. R. Garcia}, {\sc J. Mashreghi}, and {\sc W. T. Ross}, Finite Blaschke Products and Their Connections, Springer, Cham, 2018.

%\bibitem{Hardy-Littlewood-MZ-1932} {\sc G . H. Hardy} and {\sc J. E. Littlewood}, Some properties of fractional integrals. II, {\it Math. Z.} {\bf 34} (1932), 403–439.

\bibitem{Hamada-JMAA-2021} {\sc H. Hamada}, Bohr phenomenon for analytic functions subordinate to starlike or convex functions,
{\it J. Math. Anal. Appl.} {\bf 499}(1) (2021),125019.

\bibitem{Hamada-Honda-Kohr-AMP-2025} {\sc H. Hamada, T. Honda,} and {\sc M. Kohr,} Bohr–Rogosinski radius for holomorphic mappings with values in higher dimensional complex Banach spaces, \textit{Anal. Math. Phys.} (2025) 15:64, https://doi.org/10.1007/s13324-025-01061-x.

\bibitem{Ismagilov-Kayumov-Ponnusamy-JMAA-2020} {\sc A. A. Ismagilov}, {\sc I. R. Kayumov}, and {\sc S. Ponnusamy}, Sharp Bohr type inequality, {\it J. Math. Anal. Appl.} {\bf 489} (2020) 124147, 10 pp.

%\bibitem{Ismagilov-Kayumova-Kayumov-Ponnusamy-JMS-2021} {\sc A. A. Ismagilov}, {\sc A. V. Kayumova}, {\sc I. R. Kayumov}, and {\sc S. Ponnusamy}, Bohr inequalities in some classes of analytic functions, {\it J. Math. Sci.} {\bf 252} (3) (2021), 360-373.

\bibitem{Janowski-APM-1973} {\sc W. Janowski}, Some extremal problems for certain families of analytic functions I, {\it Ann. Polon. Math.} {\bf 28} (1973),297–326. doi: 10.4064/ap-28-3-297-326

\bibitem{Janowski-BAPSSSMAP-1973} {\sc W. Janowski}, Some extremal problems for certain families of analytic functions II, {\it Bull. Acad. Polon. Sci. Ser. Sci. Math. Astron. Phys.} {\bf 21} (1973),17–25.

\bibitem{Jenkins-MMJ-1962} {\sc J. A. Jenkins}, On a conjecture of Goodman concerning meromorphic univalent functions, {\it Mich. Math. J.} {\bf 9} (1962) 25–27.

\bibitem{Kalaj-MZ-2008} {\sc D. Kalaj}, Quasiconformal harmonic mapping between Jordan domains, {\it Math. Z.} {\bf 260}(2) (2008), 237–252.

\bibitem{Kaplan-MMJ-1952} {\sc W. Kaplan}, Close-to-convex schlicht functions, {\it Mich. Math. J.} {\bf 1} (1952), 169–185. 

\bibitem{Kayumov-Khammatova-Ponnusamy-JMAA-2021} {\sc I. R. Kayumov}, {\sc D. M. Khammatova}, and {\sc S. Ponnusamy}, Bohr–Rogosinski phenomenon for analytic functions and Ces$\acute{a}$ro operators, {\it J. Math. Anal. Appl.} {\bf 496} (2021) 124824, 17 pp.

\bibitem{Kayumov-Ponnusamy-CMFT-2017} {\sc I. R. Kayumov} and {\sc S. Ponnusamy}, Bohr inequality for odd analytic functions, {\it Comput. Methods	Funct. Theory} {\bf 17} (2017), 679–688.

%\bibitem{Kayumov-Ponnusamy-JMAA-2018} {\sc I. R. Kayumov} and {\sc S. Ponnusamy}, Bohr's inequalities for the analytic functions with lacunary series and harmonic functions, {\it J. Math. Anal. Appl.} {\bf 465} (2018), 857-871.

%\bibitem{Kayumov-Ponnusamy-CRMASP-2018} {\sc I. R. Kayumov} and {\sc S. Ponnusamy}, Improved version of Bohr's inequality, {C. R. Math. Acad. Sci. Paris} {\bf 356} (3) (2018), 272-277.

\bibitem{Kayumov-Ponnusamy-AASFM-2019} {\sc I. R. Kayumov} and {\sc S. Ponnusamy}, On a powered Bohr inequality, {\it Ann. Acad. Sci. Fenn. Math.} {\bf 44} (2019), 301–310.

\bibitem{Kayumov-Ponnusamy} {\sc I. R. Kayumov} and {\sc S. Ponnusamy}, Bohr-Rogosinski radius for analytic functions, preprint,
https://doi.org/10.48550/arXiv.1708.05585.

%\bibitem{Kayumov-Khammatova-Ponnusamy-CRMASP-2020} {\sc I. R. Kayumov}, {\sc D. M. Khammatova}, and {\sc S. Ponnusamy}, On the Bohr inequality for the Ces$\acute{a}$ro operator, {\it C. R. Math. Acad. Sci. Paris} {\bf 358} (2020), 615–620.

\bibitem{Kayumov-Ponnusamy-Shakirov-MN-2018}  {\sc I. R. Kayumov}, {\sc S. Ponnusamy}, and {\sc N Shakirov}, Bohr radius for locally univalent harmonic mappings, {\it Math. Nachr.} {\bf 291} (11-12) (2018), 1757-1768.

%\bibitem{Kumar-Sahoo-MJM-2021} {\sc S. Kumar} and {\sc S. K. Sahoo}, Bohr inequalities for certain integral operators, {\it Mediterr. J. Math.} {\bf 18} (2021) Article No. 268.

\bibitem{Lewy-BAMS-1936} {\sc H. Lewy}, On the non-vanishing of the Jacobian in certain one-to-one mappings, {\it Bull. Am. Math. Soc.} {\bf 42} (1936), 689–692.

\bibitem{Liu-Ponnusamy-BMMSS-2019} {\sc Z. H. Liu} and {\sc S. Ponnusamy}, Bohr radius for subordination and $K$-quasiconformal harmonic
mappings, {\it Bull. Malays. Math. Sci. Soc.} {\bf 42} (2019), 2151–2168.


%\bibitem{Liu-Liu-Ponnusamy-BDSM-2021} {\sc G. Liu}, {\sc Z. Liu}, and {\sc S. Ponnusamy}, Refined Bohr inequality for bounded analytic functions, {\it Bull. Sci. Math.} {\bf 173} (2021) 103054.

\bibitem{Liu-Ponnusami-Wang-2020} {\sc M. S. Liu}, {\sc S. Ponnusamy}, and {\sc J. Wang}, Bohr’s phenomenon for the classes of Quasisubordination
and K-quasiregular harmonic mappings, {\it Rev. Real Acad. Cienc. Exactas Fis. Nat.- A: Mat.} {\bf 114} (2020), 115.

\bibitem{Long-Wang-Wu-2022} {\sc B. Long}, {\sc Q. Wang}, and {\sc L. Wu}, some Bohr-type inequalities with one parameter for bounded analytic functions, {\it Rev. Real Acad. Cienc. Exactas Fis. Nat.- A: Mat.} {\bf 116}(2) (2022), 61.

\bibitem{Ma-Minda-1992} {\sc W. C. Ma} and {\sc D. Minda}, A unified treatment of some special classes of univalent functions. In: Proceedings
of the Conference on Complex Analysis (Tianjin, 1992), Conf. Proc. Lecture Notes Anal. I, Int. Press, Cambridge, p. 157–169.

\bibitem{Martio-AASFAI-1968} {\sc O. Martio}, On harmonic quasiconformal mappings, {\it Ann. Acad. Sci. Fenn. A. I.} {\bf 425} (1968), 3–10.

 %\bibitem{Miller-Mocanu-2000} {\sc S. S. Miller} and {\sc P. T. Mocanu}, Differential subordinations (2000) (New York: Marcel Dekker Inc.)

\bibitem{Paulsen-Popescu-Singh-PLMS-2002} {\sc V. I. Paulsen}, {\sc G. Popescu}, and {\sc D. Singh}, On Bohr's inequality, {\it Proc. Lond. Math. Soc.} {\bf 85} (2) (2002), 493-512.

\bibitem{Paulsen-Singh-2022} {\sc V. I. Paulsen} and {\sc D. Singh}, A simple proof of Bohr’s inequality, available from $arXiv:2201.10251v1$.

\bibitem{Ponnusamy-Wirths-CMFT-2020} {\sc S. Ponnusamy} and {K. -J. Wirths}, Bohr type inequalities for for functions with a multiple zero at the origin, {\it Comput. Methods Funct. Theory} {\bf 20} (2020), 559-570.

\bibitem{Robertson-AM-1936} {\sc M. S. Robertson}, On the theory of univalent functions, {\it Ann. Math.} {\bf 37} (1936),374–408. doi:
10.2307/1968451

\bibitem{Rogosinski-MJ-1923} {\sc  W. Rogosinski}, Uber Bildschranken bei Potenzreihen und ihren Abschnitten, {\it Math. Z.} {\bf 17} (1923), 260-276.

\bibitem{Wirths-SMJ-2003} {\sc K. -J. Wirths}, The Koebe domain for concave univalent functions, {\it Serdica Math. J.} {\bf 29} (2003), 355–360

\bibitem{Wirths-SMJ-2006} {\sc K. -J. Wirths}, On the residuum of concave univalent functions, {\it  Serdica Math. J.} {\bf 32} (2006), 209–214.

%\bibitem{Stempak-PRSE-1994} {\sc K. Stempak}, Ces$\acute{a}$ro averaging operators, {\it Proc.Roy.Soc.Edinburgh Sect.A} {\bf 124} (1994), 121–126.



\end{thebibliography}
\end{document}